\documentclass[12pt]{paper}
\usepackage[latin9]{inputenc}
\usepackage[a4paper]{geometry}
\usepackage{color}
\usepackage{amsmath}
\usepackage{amssymb}
\usepackage{graphicx}
\PassOptionsToPackage{normalem}{ulem}
\usepackage{ulem}

\makeatletter
\usepackage[usenames,dvipsnames,svgnames,table]{xcolor}
\usepackage{epstopdf}%\DeclareGraphicsRule{.tif}{png}{.png}{`convert #1 `dirname #1`/`basename #1 .tif`.png}

\renewcommand{\thefootnote}{\fnsymbol{footnote}}

\title{Stability of plane wave solutions in complex Ginzburg-Landau equation
with delayed feedback\thanks{We acknowledge the German Research Foundation (DFG) for financial support in the framework of the Collaborative Research Center SFB 910, and Erasmus Mundus TRIPLEI2011332 mobility program grant. A.~G. Vladimirov acknowledges the support from the SFI E.T.S. Walton award and FP7 ITN PROPHET, grant 264687.}}

\author{D. Puzyrev\footnotemark[2] \and S. Yanchuk\footnotemark[2]
        \and A.~G. Vladimirov\footnotemark[3] \footnotemark[4] \and S.~V. Gurevich\footnotemark[5]}

\begin{document}

\maketitle

\footnotetext[2]{Institute of Mathematics, Humboldt University of Berlin, Unter den Linden 6, D-10099 Berlin, Germany ({\tt puzyrev@math.hu-berlin.de}).}
\footnotetext[3]{Weierstrass Institute for Applied Analysis and Stochastics, Mohrenstrasse 39, D-10117 Berlin, Germany}
\footnotetext[4]{Cork Institute of Technology, Rossa Ave, Bishopstown, Cork, Ireland}
\footnotetext[5]{Institute for Theoretical Physics, University of M\"unster, Wilhelm-Klemm-Str. 9, D-48149 M\"unster, Germany}

\renewcommand{\thefootnote}{\arabic{footnote}}

\begin{abstract}
We perform bifurcation analysis of plane wave solutions in one-dimensional
\textcolor{black}{complex} cubic-quintic Ginzburg-Landau equation
with delayed feedback. Our study reveals how multistability and snaking
behavior of plane waves emerge as time delay is introduced. For intermediate
values of the delay, bifurcation diagrams are obtained by a combination
of analytical and numerical methods. For large delays, using an asymptotic
approach we classify plane wave solutions into strongly unstable,
weakly unstable, and stable. The results of analytical bifurcation
analysis are in agreement with those obtained by direct numerical
integration of the model equation. \end{abstract}
\begin{keywords}
delay differential equation, stability theory, large delay, periodic
solutions.
\end{keywords}
%\begin{AMS}
%34K08, 34K20, 35R10
%\end{AMS}

\section{Introduction}

The complex Ginzburg-Landau equation (CGLE) plays an important role
in modeling of various natural phenomena including nonlinear optical
waves, second-order phase transitions, Rayleigh-Bénard convection,
and superconductivity \cite{RevModPhys.65.851,Aranson-Kramer2002}.
It is an amplitude equation describing the onset of instability near
an Andronov-Hopf bifurcation in spatially extended dynamical systems
\cite{newell74}. In nonlinear optics, equations of CGLE type are
widely used to describe such phenomena as mode-locking in lasers \cite{soto96,akhmediev98,Grelu2004},
short pulse propagation in optical transmission lines \cite{matsumoto95},
dynamics of multimode lasers, and transverse pattern formation in
nonlinear optical media \cite{VFKKR99,Staliunas2007}. While classical
cubic CGLE describes a supercritical bifurcation, in the case of subcritical
instability this equation is augmented with the fifth-order terms
to allow the existence of stable pulselike solutions \cite{Thual88,kivshar2003optical}.

Here we focus on one-dimensional delayed cubic-quintic CGLE for the
complex amplitude $A(x,t)$ 
\begin{equation}
\partial_{t}A=\left(\beta+\frac{i}{2}\right)\partial_{xx}A+\delta A+\left(\epsilon+i\right)|A|^{2}A+\left(\mu+i\nu\right)|A|^{4}A+\eta e^{i\varphi}A(x,t-\tau).\label{CGLE}
\end{equation}
The parameter $\beta>0$ is the diffusion coefficient, second-order
dispersion (diffraction) is scaled to $1/2$, and $\delta$ describes
the linear loss or gain. Parameters $\epsilon$, $\mu$, and $\nu$
determine the shape of the nonlinearity. In particular, there are
two qualitatively different cases: $\epsilon>0$ and $\epsilon<0$
corresponding to destabilizing and stabilizing role of the cubic nonlinearity.
Equivalently, this leads to either subcritical or supercritical destabilization
mechanisms for the homogeneous steady state $A=0$. In this work, we
take the values $\varepsilon<0$, $\mu=\nu=0$ for the supercritical
case (cubic CGLE), and $\varepsilon>0$, $\mu<0$ for the subcritical
case. Parameters $\eta$ and $\varphi$ determine the feedback rate
and phase, respectively, whereas $\tau$ is the delay time. Model
equation \eqref{CGLE} can describe, for instance, a broad area optical
system with delayed optical feedback \cite{tlidi09,sveta13}. Notice
that in the absence of delayed feedback term, $\eta=0$, Eq. \eqref{CGLE}
becomes the classical cubic-quintic CGLE \cite{RevModPhys.65.851,hohenberg1992Fronts}.

Although CGLE possesses a variety of different solutions \cite{Akh-Soto1,Soto-Crespo2,Aranson-Kramer2002,Mielke2002},
in this work we focus on the simplest plane waves of the form $A=a_{0}e^{iqx+i\omega t}$.
The stability of the trivial homogeneous solution $A=0$ is studied
as well. Criteria of the stability of plane wave solutions in the
quintic CGLE without delay were briefly described in \cite{hohenberg1992Fronts}.
The stabilization of plane waves in one-dimensional and two-dimensional
cubic CGLE by a combination of spatially shifted and temporally delayed
non-invasive feedback was investigated in \cite{Montgomery04feedbackcontrol,Postlethwaite2007}
for the case when delay time and space shift are in the resonance
with plane wave spatial and temporal wavenumbers. In this paper, we
study cubic-quintic CGLE (as well as cubic CGLE as a special case)
with arbitrary delay time and phase of the feedback, including the
long delay limit case. For small delay times, there appears a single
plane wave for every allowed spatial wavenumber $q$. The local stability
of such solutions can be studied by calculating a dispersion relation
for a given plane wave solution \cite{hohenberg1992Fronts}. As time
delay becomes comparable or longer than the characteristic time scale
in the system the number of plane waves for each admissible wavenumber
grows linearly with $\tau$. Moreover, the stability of each plane
wave is no longer determined by a single classical dispersion relation,
but a set of dispersion relations, which correspond to various ``delay-induced
modes''. This set of dispersion relations is conveniently described
using the methods from \cite{Wolfrum2010,Yanchuk2010a,Lichtner2011,Sieber2013,Giacomelli1996}
by adding an additional dimension to the dispersion relation and studying
a so called ``hybrid dispersion relation'', which is a function
of two arguments. Performing such a stability analysis, we identify
a large set of emerging asymptotically stable plane waves. Another
contribution of this work is that we present a way how one can conveniently
describe and visualize the whole set of plane waves and their stability
in system \eqref{CGLE}. As an interesting observation, we find that
the branches of plane waves exhibit a snaking behavior as the linear
gain parameter $\delta$ is changed.

This paper is organized as follows: in Section 2 we start with the
stability analysis of the homogeneous solution $A=0$. In addition,
in this section we introduce some important ideas that will be used
in a technically more elaborated way in Section 3, where the existence
and stability of the plane wave solutions $A=a_{0}e^{iqx+i\omega t}$
is studied. Finally, the conclusions are given in Section 4.

\section{Stability analysis of the homogeneous solution}

\subsection{The case without delayed feedback}

In this section we start with the stability analysis of the trivial
homogeneous steady state solution $A(x,t)=0$. Let us briefly recall
the case when the feedback is absent, i.e. $\eta=0$ \cite{hohenberg1992Fronts,Montgomery04feedbackcontrol}.
By substituting the perturbations of the form $e^{iqx+\lambda t}$
in the linearized equation, we obtain the characteristic equation
\begin{equation}
\chi(\lambda)=\lambda-\delta+\left(\beta+\frac{i}{2}\right)q^{2}=0.\label{CP1}
\end{equation}
Here $q$ is a spatial wavenumber of the perturbation and $\lambda$
determines the growth rate. The corresponding dispersion relation
reads 
\[
\lambda(q)=\delta-\left(\beta+\frac{i}{2}\right)q^{2}.
\]
When all the eigenvalues $\lambda(q)$ have negative real parts, the
homogeneous solution is asymptotically stable. Since $\beta$ is positive,
we conclude that the trivial solution is unstable for $\delta>0$
and stable if $\delta<0$. The most unstable wavenumber is $q=0$.

\subsection{Case of delay $\tau$}

For nonzero feedback rate the characteristic equation for the homogeneous
solution becomes
\begin{equation}
\chi_{1}(\lambda)=\lambda-\delta-\eta e^{i\varphi}e^{-\lambda\tau}+\left(\beta+\frac{i}{2}\right)q^{2}=0.\label{CP2}
\end{equation}

\begin{figure}[h]
\centering \includegraphics[width=0.2\linewidth]{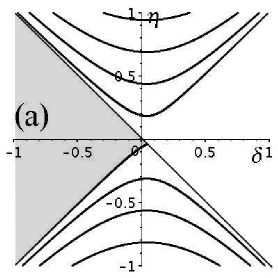}\hskip
0.5 cm\includegraphics[width=0.2\linewidth]{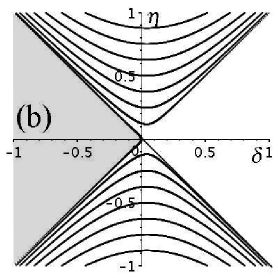}\hskip
0.5 cm\includegraphics[width=0.2\linewidth]{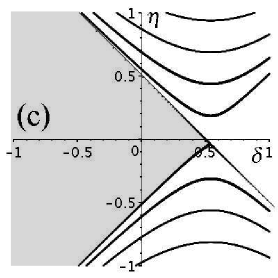}\hskip
0.5 cm\includegraphics[width=0.2\linewidth]{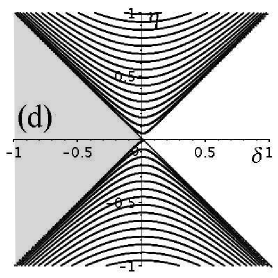}

\caption{Andronov-Hopf bifurcation curves for the trivial homogeneous solution
$A=0$ of Eq. \eqref{CGLE} in $(\delta,\eta)$ plane. The curves
are given in a parametric form \eqref{param2}. Stability region of
the solution $A=0$ where the real part of the eigenvalue corresponding
to the most unstable mode $q=0$ is negative is shown in gray. Parameter
values: (a) $\tau=25,\, q=0,\,\varphi=0,$ (b) $\tau=25,\, q=0,\,\varphi=\pi/2,$
(c) $\tau=25,\, q=1,\,\varphi=0,$ (d) $\tau=100,\, q=0,\,\varphi=0$.
In all figures $\beta=0.5$.\label{fig:bif2}}
\end{figure}

The Andronov-Hopf bifurcation curves in the parameter plane $(\eta,\delta)$
can be found in a parametric form. To this end we substitute $\lambda=i\omega$
into characteristic equation \eqref{CP2} 
\begin{equation}
i\omega-\delta-\eta e^{i\varphi}e^{-i\omega\tau}+\left(\beta+\frac{i}{2}\right)q^{2}=0.\label{CP2-1}
\end{equation}
By separating real and imaginary parts of equation \eqref{CP2-1},
we obtain two relations 
\begin{equation}
\eta(\omega)=\frac{\frac{q^{2}}{2}+\omega}{\sin(\varphi-\omega\tau)},\quad\delta(\omega)=\beta q^{2}-\frac{\cos(\varphi-\omega\tau)\left(\frac{q^{2}}{2}+\omega\right)}{\sin(\varphi-\omega\tau)}\label{param2}
\end{equation}
defining the Andronov-Hopf bifurcation curves with the imaginary part
of the critical eigenvalue $\omega$ being a free parameter. Figure
\ref{fig:bif2} shows these bifurcation curves in the plane of two
parameters, linear gain $\delta$ and feedback rate $\eta$. The stability
region of the trivial solution where real parts of all the eigenvalues
$\lambda(q)$ are negative is shown in gray. A destabilization with
respect to a given wavenumber $q$ takes place when crossing the boundary
of this region from inside. Figures \ref{fig:bif2}(a,b) present the
Andronov-Hopf bifurcation curves for the wavevector $q=0$ at different
values of the feedback phase, see also \cite{ReddyOsc,Montgomery04feedbackcontrol}.
For non-zero values of $q$, the instability region shifts to higher
values of $\delta$, as shown in Fig. \ref{fig:bif2}(c). Therefore,
the destabilization of the trivial homogeneous steady state first
occurs at the most unstable wavenumber $q=0.$ For larger delay times
the set of bifurcation curves becomes more dense, as it is seen from
Fig. \ref{fig:bif2}(d).

\subsection{The case of large delay\label{sub:homlong}}

Delay time $\tau$ for optical systems can exceed the internal time
scales of the system by several orders of magnitude \cite{Soriano2013,Fiedler2008,Larger2013,Giacomelli2012}.
To study this case, we use the asymptotic technique similar to that
used in \cite{Wolfrum2010,Yanchuk2010a,Lichtner2011,Sieber2013}.
According to these results, there are two types of delay-induced instabilities:
strong instability and weak instability. Strong instability appears
when there exists an eigenvalue (or Lyapunov exponent, for the periodic
or chaotic case), which tends to some constant value $\lambda(\tau)\to\lambda_{0}$
with $\Re[\lambda_{0}]>0$ as time delay increases. In this case,
the contribution of the term $e^{-\lambda\tau}$ in the characteristic
equation \eqref{CP2} can be neglected and we arrive to the following
condition for the strong instability of the mode with the wavenumber
$q$ 
\[
\lambda_{0}=\delta+\left(\beta+\frac{i}{2}\right)q^{2}>0.
\]
This instability condition coincides, in fact, with that for the feedback-free
case. The most unstable mode $q=0$ gives the condition for the strong
instability of the solution $A=0$: 
\begin{equation}
\delta>0.\label{StrInst}
\end{equation}

Another type of instability, the weak one, can be described by a pseudo-continuous
spectrum of eigenvalues, which scales as 
\begin{equation}
\lambda=\frac{\gamma(\xi)}{\tau}+i\xi,\label{PCS}
\end{equation}
in the limit $\tau\to\infty$. More specifically, it was proved in
\cite{Lichtner2011,Sieber2013} that this spectrum is converging to
a set of continuous curves \eqref{PCS} parametrized by $\xi$. Moreover,
the leading terms of the real parts $\gamma(\xi)$ can be found explicitly
by substituting \eqref{PCS} into the characteristic equation \eqref{CP2}
and neglecting small terms of order $1/\tau$. In our case, we obtain
the equation

\begin{equation}
i\xi-\delta-\eta e^{i\varphi}e^{-\gamma}e^{-i\xi\tau}+\left(\beta+\frac{i}{2}\right)q^{2}=0,\label{dPCS}
\end{equation}
which is solved explicitly with respect to $\gamma$: 
\begin{equation}
\gamma(\xi,q)=-\frac{1}{2}\ln\frac{\left(\delta-\beta q^{2}\right)^{2}+\left(\xi+\frac{1}{2}q^{2}\right)^{2}}{\eta^{2}}.\label{gamma}
\end{equation}
In Eq. \eqref{gamma}, we write an additional argument $q$, which
indicates the dependence of the real part (rescaled by $1/\tau$)
of the eigenvalues on the wavenumber. The function $\gamma$ of two
arguments extends naturally the dispersion relation, which is used
\cite{Cross2009} for partial differential equations without delayed
feedback. Indeed, for a fixed $\xi$, the relation \eqref{gamma}
determines the stability of the homogeneous state with respect to
the perturbations with the wavenumber $q$, i.e. it is the dispersion
relation. The new variable $\xi$ stands for the delay induced modes,
which appear additionally due to the delay. The homogeneous solution
is locally asymptotically stable when $\gamma(\xi,q)<0$.

\begin{figure}[h]
\centering \includegraphics[width=0.3\linewidth]{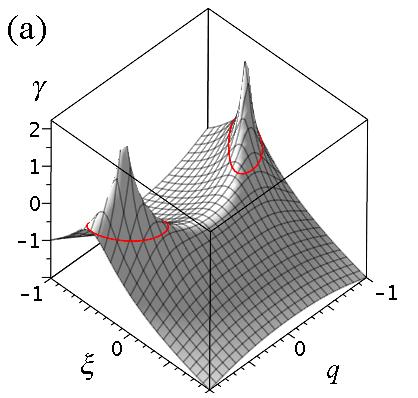}
\includegraphics[width=0.3\linewidth]{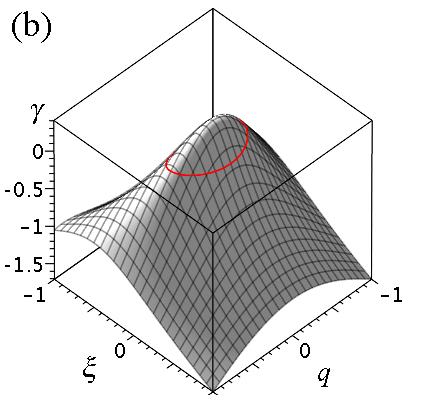}\includegraphics[width=0.3\linewidth]{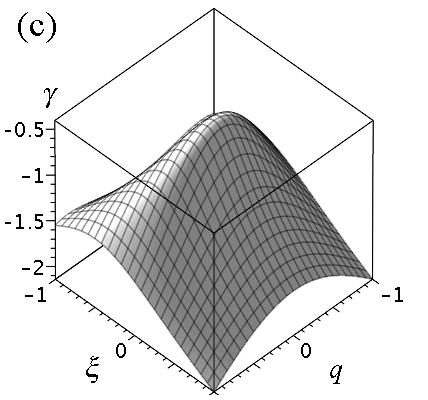}
\caption{Dispersion relation $\gamma(\xi,q)$, given by Eq. \eqref{gamma}
for the homogeneous state $A=0$ at large $\tau$. Here $\xi$ and
$q$ stand for the delay induced and the spatial wavenumber of the
perturbations. If $\gamma<0$ for all $\xi$ and $q$ and no strong
instability occurs, the homogeneous state is stable. Parameter values:
(a) unstable case with two regions of unstable wavenumbers: $\delta=0.3,\,\eta=0.2$,
(b) unstable case with one unstable region: $\delta=-0.2,\,\eta=0.3$,
(c) stable case: $\delta=-0.3,\,\eta=0.2$. Red curve shows the level
lines $\gamma=0$. In all figures $\beta=0.5$.\label{fig:sol} }
\end{figure}

\begin{figure}[h]
\centering \includegraphics[width=0.4\linewidth]{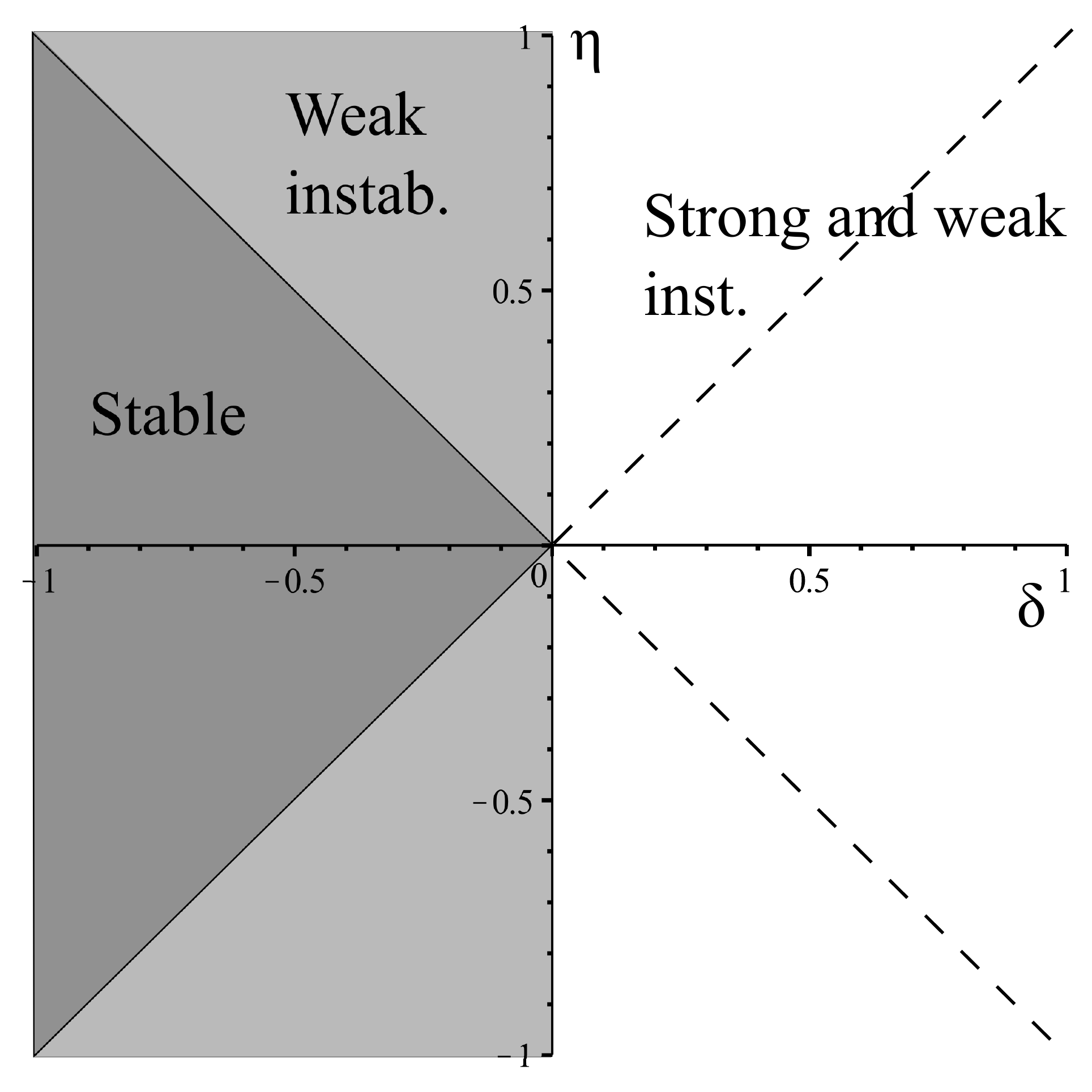}
\caption{Bifurcation diagram of the homogeneous solution $A=0$ in the parameter
plane $(\delta,\eta)$ for large delay times. The diagram shows the
region of stability (gray), weak and strong instability (light gray
and white, respectively). \label{fig:bif1} }
\end{figure}

Figure \ref{fig:sol} shows the surfaces of $\gamma(\xi,q)$ calculated
for different parameters. The red curve shows the level line $\gamma=0$
given by the relation 
\begin{equation}
\left(\delta-\beta q^{2}\right)^{2}+\left(\frac{q^{2}}{2}+\xi\right)^{2}=\eta^{2}.\label{edl1}
\end{equation}
Equation \eqref{edl1} defines an ellipse in $\xi$ and $q^{2}$ coordinates.
The trivial solution $A=0$ is unstable if at least part of the ellipse
lies in the half-plane $q^{2}>0$. Simple calculations show that,
for $\beta>0$, this is equivalent to the condition 
\begin{equation}
\delta>-|\eta|.\label{cond2}
\end{equation}
Hence, the inequality \eqref{cond2} gives the weak instability condition
for the solution $A=0$, see Fig. \ref{fig:bif1}. The corresponding
critical wave numbers $\xi_{c}$ and $q_{c}$ correspond to the maximum
of the quantity $\gamma(\xi,q)$. For $\delta\le0$, we have $q_{c}=0$
and $\xi_{c}=0.$ Under the condition 
\begin{equation}
\delta>|\eta|\label{cond3}
\end{equation}
there are two separated regions of unstable wavenumbers $\xi$ and
$q$, for which $\gamma(\xi,q)>0$, see Fig. \ref{fig:sol}(a). These
regions are symmetric with respect to the $\xi$-axis, $q=0$. Otherwise,
when only first of the two inequalities, \eqref{cond2} and \eqref{cond3},
is satisfied, there is a single symmetric region of unstable wavenumbers,
see Fig.~\ref{fig:sol}(b). The boundary defined by the inequality
\eqref{cond3} is shown in Fig.~\ref{fig:bif1} by a dashed line.
Figure \ref{fig:sol}(c) illustrates the case when the eigenvalue
spectrum is located in the left half-plane of the complex plane, $\gamma(\xi,q)<0$,
and the homogeneous state $A=0$ is stable.

The complete bifurcation diagram for the homogeneous state in the
case of long delay is summarized in Fig.~\ref{fig:bif1}, where the
regions of strong and weak instability are shown. It is instructive
to compare this bifurcation diagram obtained in the limit $\tau\to0$
with the exact bifurcation curves for different values of $\tau$
shown at Figs. \ref{fig:bif2}(a,b,d). 

\section{Plane wave solutions}

In this section we study existence and stability of plane wave solutions
$A=a_{0}e^{iqx+i\omega t}$ of Eq. \eqref{CGLE}. As for the homogeneous
solution, we start with the case without delayed feedback in Sec.~\ref{sub:Stability-of-pw},
and then consider the case with finite delay in Sec.~\ref{sub:Case-of-delay}.
Further, in Sec.~\ref{sub:The-case-of-large} we study analytically
the limit of a large delay time $\tau$, which allows for a deeper
insight into the stability properties of plane wave solutions as compared
to the arbitrary $\tau$ case. Finally, in Sec. \ref{sub:Numerical-simulations-of}
the results of numerical simulations are presented.

\subsection{The case without delayed feedback \label{sub:Stability-of-pw}}

Substituting $A=a_{0}e^{iqx+i\omega t}$ into Eq. \eqref{CGLE}, we
obtain the relation between the unknown amplitude $a_{0}$, wavenumber
$q$, and frequency $\omega$ of the plane wave solution: 
\begin{equation}
i\omega=-\left(\beta+\frac{i}{2}\right)q^{2}+\delta+\left(\epsilon+i\right)a_{0}^{2}+\left(\mu+i\nu\right)a_{0}^{4}.\label{eq:pwaves1}
\end{equation}
Due to the symmetry property of the CGLE $A\to-A$, this equation
is symmetric under the reflection $a_{0}\to-a_{0}$. Hence, we restrict
our analysis to the case $a_{0}\geq0$. The real and imaginary parts
of Eq. \eqref{eq:pwaves1} give the expressions for the amplitude
$a_{0}^{2}(q)$ and the frequency $\omega(q)$ at a given wavenumber
$q$: 
\begin{equation}
a_{0}^{2}(q)=\frac{-\epsilon\pm\sqrt{\epsilon^{2}-4\mu(\delta-\beta q^{2})}}{2\mu},\quad\omega(q)=-\frac{q^{2}}{2}+a_{0}^{2}(q)+\nu a_{0}^{4}(q).\label{eq:wsol1}
\end{equation}
In particular, for the cubic CGLE (supercritical case with $\mu=\nu=0$)
we obtain 
\begin{equation}
a_{0}^{2}(q)=\sqrt{\frac{\beta q^{2}-\delta}{\epsilon}},\quad\omega(q)=-\frac{q^{2}}{2}+\frac{\beta q^{2}-\delta}{\epsilon}.\label{eq:wsol1cubic}
\end{equation}

\begin{figure}[h]
\centering \includegraphics[width=0.3\linewidth]{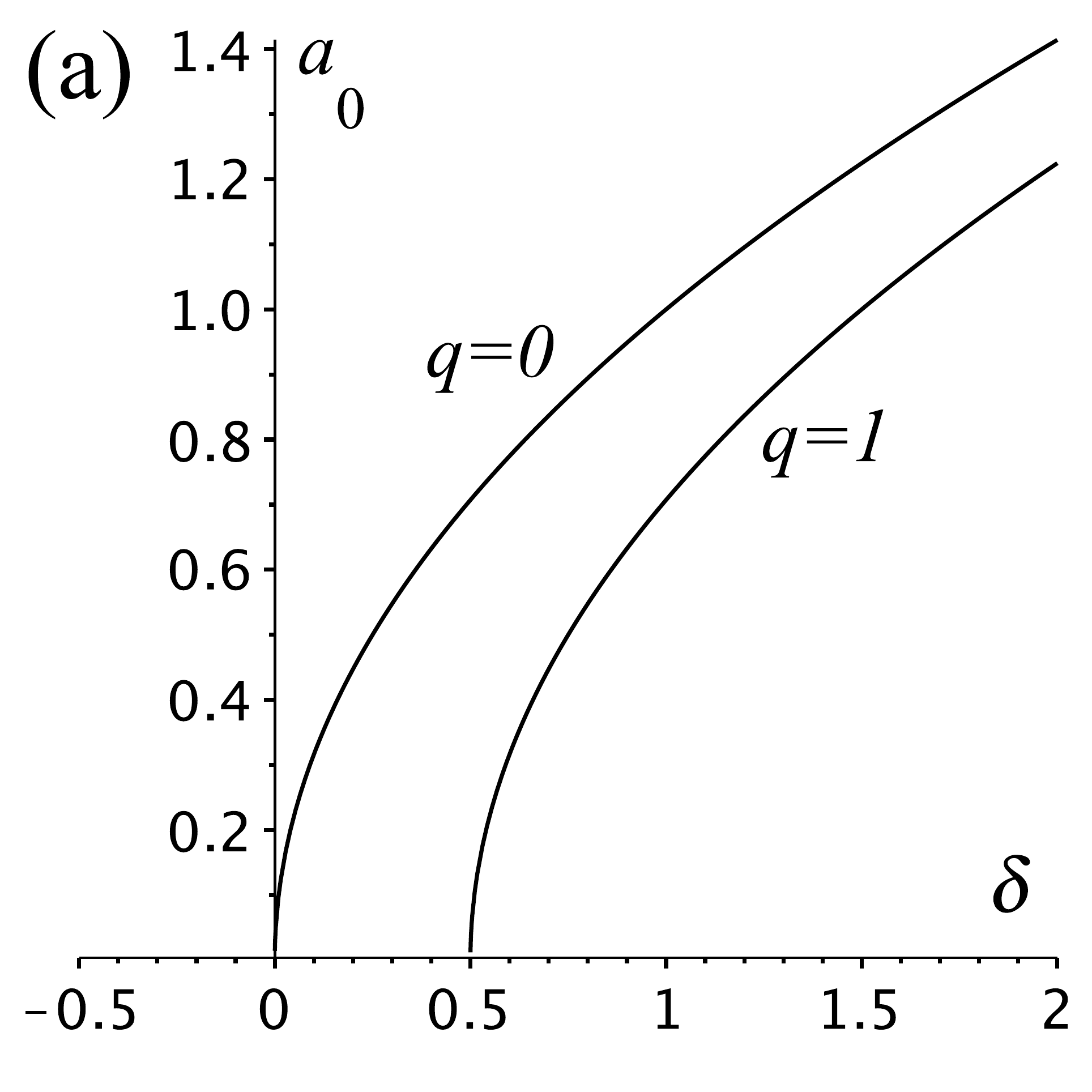}\hskip
1.5 cm\includegraphics[width=0.3\linewidth]{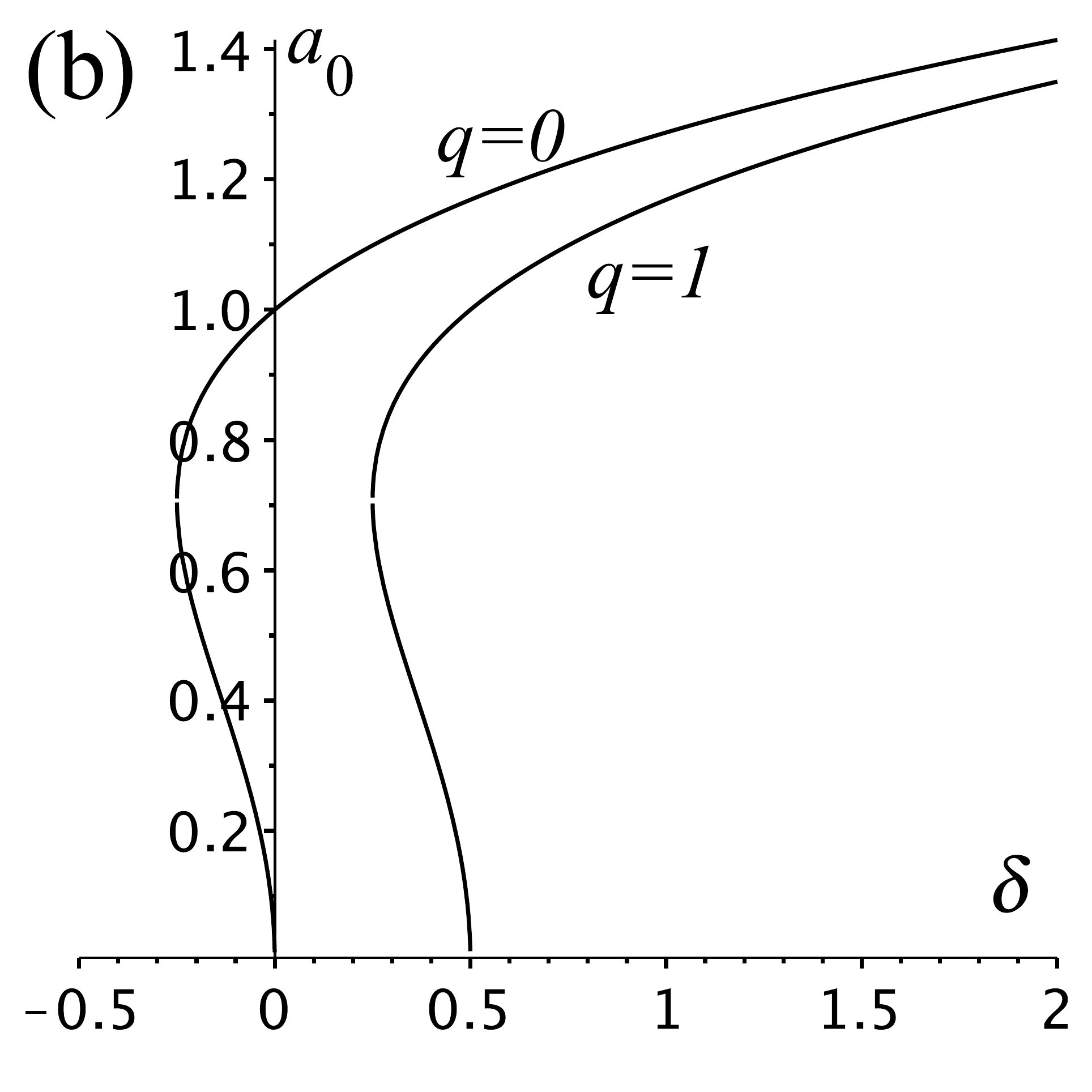}
\caption{Amplitude of plane wave solutions $a_{0}(\delta)$ versus parameter
$\delta$ in the absence of the delay, given by Eqs. \eqref{eq:wsol1cubic}
and \eqref{eq:wsol1}, for different spatial wavenumbers $q$. (a)
Supercritical case, $\epsilon=-1,$ $\mu=\nu=0$. (b) Subcritical
case, $\epsilon=1$, $\mu=-1$, $\nu=-0.1$. In all plots $\beta=0.5.$\label{fig:wsol1} }
\end{figure}

Figures \ref{fig:wsol1} (a) and (b) show the amplitude of the plane
wave $a_{0}(\delta)$ as a function of the gain parameter $\delta$
for the supercritical and subcritical case, respectively. The branch
of plane wave solutions with a given $q$ emerges from the homogeneous
state via Hopf bifurcation at $\delta=\beta q^{2}$. For the parameter
values of Fig. \ref{fig:wsol1}(b) corresponding to the quintic CGLE,
it bifurcates subcritically from $A=0$ and undergoes a fold bifurcation
at $\delta=\frac{\epsilon^{2}}{4\mu}+\beta q^{2}.$

\begin{figure}[h]
\centering \includegraphics[width=0.6\linewidth]{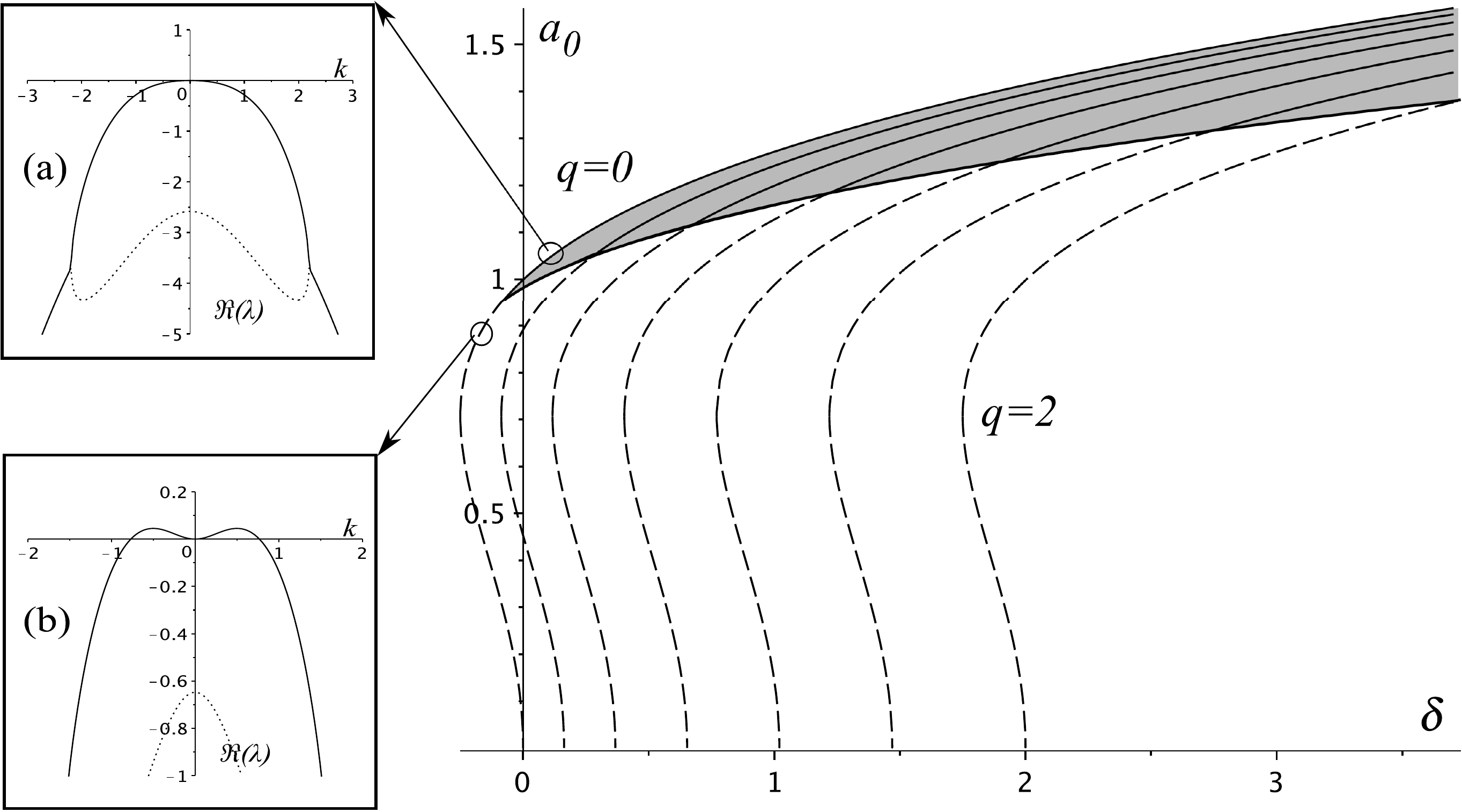}
\caption{Amplitude of plane wave solutions $a_{0}(\delta)$ in quintic CGLE
\eqref{CGLE} in the absence of the delay for different spatial wavenumbers
$q$. Stability domain is shown in gray. Insets show the dependence
of real part of the growth rate $\lambda(k)$ for selected plane waves.
(a) $\delta=0.1,\, q=0,\, a_{0}\simeq1.04$, stable case, (b) $\delta=-0.2,\, q=0,\, a_{0}\simeq0.85$,
unstable case. Other parameters: $\beta=0.5,$ $\epsilon=1,$ $\mu=-1$,
and $\nu=-0.1$.\label{fig:Ekh3}}
\end{figure}

Note that plane wave solutions of the CGLE without delay in supercritical
case show the classical Benjamin-Feir scenario \cite{hohenberg1992Fronts,Aranson-Kramer2002}.
Let us shortly discuss the stability of plane waves in the subcritical
case of the quintic CGLE. Although the main ideas of this analysis
are known from e.g. \cite{hohenberg1992Fronts}, there are still some
details, which are not explained in \cite{hohenberg1992Fronts}, but
are important for our further analysis. The perturbed plane wave solutions
have the form 
\begin{equation}
A(x,t)=(a_{0}+a_{p})e^{iqx+i\omega t},\label{eq:pw-ansatz}
\end{equation}
where 
\begin{equation}
a_{p}=a_{+}e^{ikx+\lambda t}+\bar{a}_{-}e^{-ikx+\overline{\lambda}t}\label{eq:ap}
\end{equation}
is a small perturbation term with a growth rate $\text{\ensuremath{\lambda}}$.
Here $\bar{\lambda}$ and $\bar{a}$ denote complex conjugation, and
$k$ stands for different perturbation modes. Substitution of \eqref{eq:pw-ansatz}
into CGLE \eqref{CGLE} with $\eta=0$ and linearization in $a_{p}$
yields 
\begin{eqnarray}
\partial_{t}a_{p}+i\omega\left(a_{0}+a_{p}\right)=\left(\beta+\frac{i}{2}\right)\left(\partial_{xx}a_{p}+2iq\partial_{x}a_{p}-q^{2}(a_{0}+a_{p})\right)+\delta\left(a_{0}+a_{p}\right)+\nonumber \\
+\left(\epsilon+i\right)\left(a_{0}^{3}+2a_{0}^{2}a_{p}+a_{0}^{2}\overline{a_{p}}\right)+\left(\mu+i\nu\right)\left(a_{0}^{5}+3a_{0}^{2}a_{p}+2a_{0}^{2}\overline{a_{p}}\right).\label{eq:pw2}
\end{eqnarray}
After substituting \eqref{eq:ap} into \eqref{eq:pw2} and using Eq.
\eqref{eq:pwaves1} we obtain an equation involving two linearly independent
functions $e^{ikx+\lambda t}$ and $e^{-ikx+\overline{\lambda}t}.$
Requiring that the coefficients at these functions are zero, we arrive
at a system of linear equations for the unknowns $a_{-}$ and $a_{+}$:
\[
M\left(\begin{array}{c}
a_{+}\\
a_{-}
\end{array}\right)=0
\]
with {\small 
\begin{equation}
M=\begin{bmatrix}\lambda+\left(\beta+\frac{i}{2}\right)\left(k^{2}+2kq\right)- & \qquad & -\left(\epsilon+i\right)a_{0}^{2}-2\left(\mu+i\nu\right)a_{0}^{4}\\
-\left(\epsilon+i\right)a_{0}^{2}-2\left(\mu+i\nu\right)a_{0}^{4}\\
\\
 &  & \lambda+\left(\beta-\frac{i}{2}\right)\left(k^{2}-2kq\right)-\\
-\left(\epsilon-i\right)a_{0}^{2}-2\left(\mu-i\nu\right)a_{0}^{4} &  & -\left(\epsilon-i\right)a_{0}^{2}-\left(\mu-i\nu\right)a_{0}^{4}
\end{bmatrix}.
\end{equation}
}Since we are looking for non-trivial solutions $\left(a_{+},a_{-}\right)$,
the characteristic equation for the perturbation growth rate $\lambda(k)$
is obtained by setting $\det M=0$: 
\begin{eqnarray}
\lambda^{2}+2\left(ikq+\beta k^{2}-\epsilon a_{0}^{2}-2\mu a_{0}^{4}\right)\lambda-2\left(\left(\nu+2\mu\beta\right)k^{2}+2\left(\mu-2\nu\beta\right)ikq\right)a_{0}^{4}+\nonumber \\
-\left(\left(1+2\epsilon\beta\right)k^{2}+2\left(\epsilon-2\beta\right)ikq\right)a_{0}^{2}+\left(4\beta^{2}+1\right)\left(k^{4}/4+k^{2}q^{2}\right)=0.\label{eq:lambda}
\end{eqnarray}
Solutions $\lambda(k)$ can now be found explicitly and the maximum
of their real parts determines the stability of plane waves. As it
is seen from Fig.~\ref{fig:Ekh3} the plane waves are stable for
larger values of $\delta$ and become modulationally unstable with
the decrease of $\delta$. This instability appears when the real
part of the derivative $\partial_{kk}\lambda(0)$ changes its sign
from negative to positive, see insets (a) and (b) in Fig. \ref{fig:Ekh3}.

The threshold of the modulational instability is given by the condition
$\Re[\partial_{kk}\lambda(0)]=0$. For small $k$ (long-wavelength
limit) the Taylor expansion of $\lambda(k)$ reads

\begin{equation}
\lambda(k)=\left(-\frac{C_{3}}{C_{1}}-q\right)ik+\left(-\frac{C_{3}^{2}}{128C_{1}^{3}}-\frac{C_{2}^{2}}{C_{1}}-\beta\right)k^{2}+\mathcal{O}(k^{3}),
\end{equation}
where $C_{1}=\epsilon a_{0}^{2}+2\mu a_{0}^{4},$ $C_{2}=16\beta^{2}q^{2}+4a_{0}^{2}+8\nu a_{0}^{4}$
and $C_{3}=64\beta^{3}q^{3}-4\beta qC_{2}.$ In particular, the stability
boundary $\partial_{kk}\lambda(0)=0$ is given by 
\begin{equation}
\frac{C_{3}^{2}}{128C_{1}^{3}}+\frac{C_{2}^{2}}{C_{1}}+\beta=0.
\end{equation}

Figure \ref{fig:Ekh3} shows the amplitude $a_{0}$ of plane wave
solutions versus $\delta$ for different wavenumbers $q$ along with
the examples of $\Re[\lambda(k)]$ for stable an unstable cases respectively.
Stable solutions are depicted by solid lines and the stable domain
is shown in gray. The first zero wavenumber plane wave solution $A=a_{0}e^{i\omega t}$
with $q=0$, $a_{0}^{2}=[-\epsilon\pm\sqrt{\epsilon^{2}-4\mu\delta}]/2\mu$,
and $\omega=a_{0}^{2}+\nu a_{0}^{4}$ is stable for 
\begin{equation}
\delta>\frac{4\beta\epsilon\left(\beta\epsilon+2\beta\mu+\nu+1\right)+2(2\beta\mu+\nu+1)}{4(2\beta\mu+\nu)^{2}}.
\end{equation}
With the further increase of $\delta$, stable plane wave solutions
with $q^{2}>0$ appear.

\subsection{Case of delay $\tau$\label{sub:Case-of-delay}}

\subsubsection*{Description of the set of plane wave solutions}

Substituting $A=a_{0}e^{iqx+i\omega t}$ into CGLE with delayed feedback
\eqref{CGLE} we obtain the equation connecting the amplitude $a_{0}$,
frequency $\omega$, and wavenumber $q$ of the plane waves 
\begin{equation}
i\omega=-\left(\beta+\frac{i}{2}\right)q^{2}+\delta+\left(\epsilon+i\right)a_{0}^{2}+\left(\mu+i\nu\right)a_{0}^{4}+\eta e^{i\varphi}e^{-i\omega\tau}.\label{eq:pwdelaymain}
\end{equation}
After separating real and imaginary parts of \eqref{eq:pwdelaymain},
we obtain 
\begin{eqnarray}
0 & = & -\beta q^{2}+\delta+\epsilon a_{0}^{2}+\mu a_{0}^{4}+\eta\cos(\omega\tau-\varphi),\label{eq:pw-delay1}\\
\omega & = & -\frac{q^{2}}{2}+a_{0}^{2}+\nu a_{0}^{4}-\eta\sin(\omega\tau-\varphi).\label{eq:pw-delay2}
\end{eqnarray}
\begin{figure}[h]
\centering \includegraphics[width=0.6\linewidth]{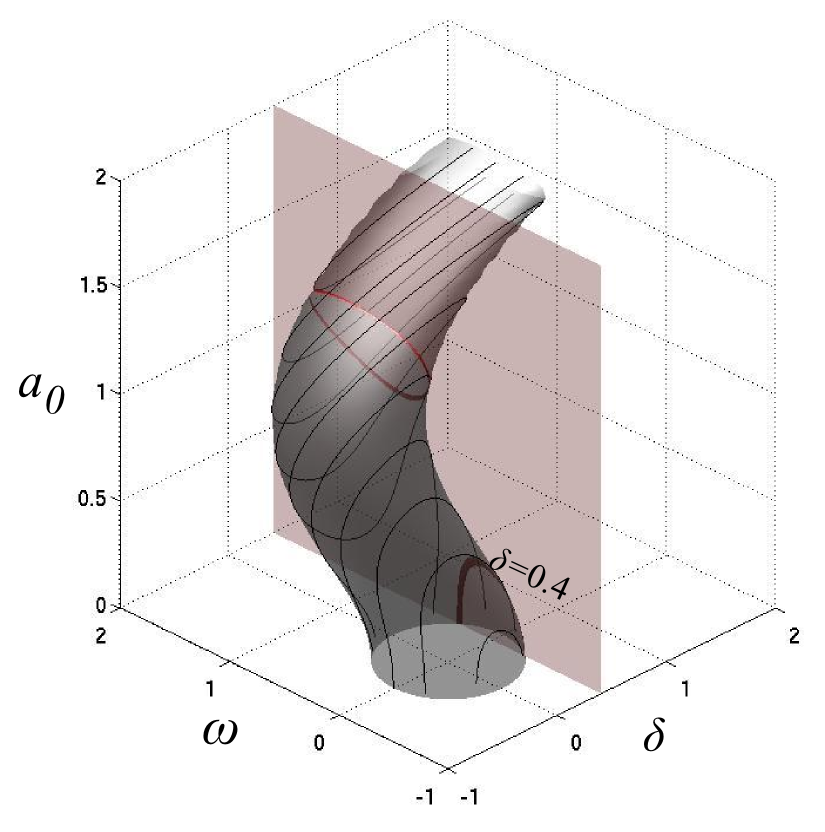}\caption{Amplitudes $a_{0}$ and frequencies $\omega$ of plane wave solutions
as a function of $\delta$, defined by Eq.~\eqref{eq:pwd-ell1}.
Black curve indicates the branch of plane wave solutions (defined
by \eqref{eq:pw-delay1}--\eqref{eq:pw-delay2}) for $\tau=25$ and
$\varphi=0$. Other parameters: $q=0,\:\eta=0.5,\:\beta=0.5,$ $\epsilon=1,$
$\mu=-1$, and $\nu=-0.1$. Red curve shows the cross-section of the
tube by the plane of the fixed parameter $\delta=0.4$.\label{fig:soltube-1}}
\end{figure}

In the coordinates $(\delta,\omega,a_{0})$, the set of solutions
of Eqs. \eqref{eq:pw-delay1}--\eqref{eq:pw-delay2} for each $q$
lies on the surface of a tube, which is implicitly defined by the
following equation 
\begin{equation}
\left(-\beta q^{2}+\delta+\epsilon a_{0}^{2}+\mu a_{0}^{4}\right)^{2}+\left(\omega+\frac{q^{2}}{2}-a_{0}^{2}-\nu a_{0}^{4}\right)^{2}=\eta^{2}\label{eq:pwd-ell1}
\end{equation}
(see Fig. \ref{fig:soltube-1}), obtained from Eqs. \eqref{eq:pw-delay1}--\eqref{eq:pw-delay2}
by exclusion of $\omega\tau-\varphi$. For fixed parameters, including
$\delta$, this is a one-dimensional set, see the cross-section in
Fig.~\ref{fig:soltube-1}.

\begin{figure}[h]
\centering \includegraphics[width=0.3\linewidth]{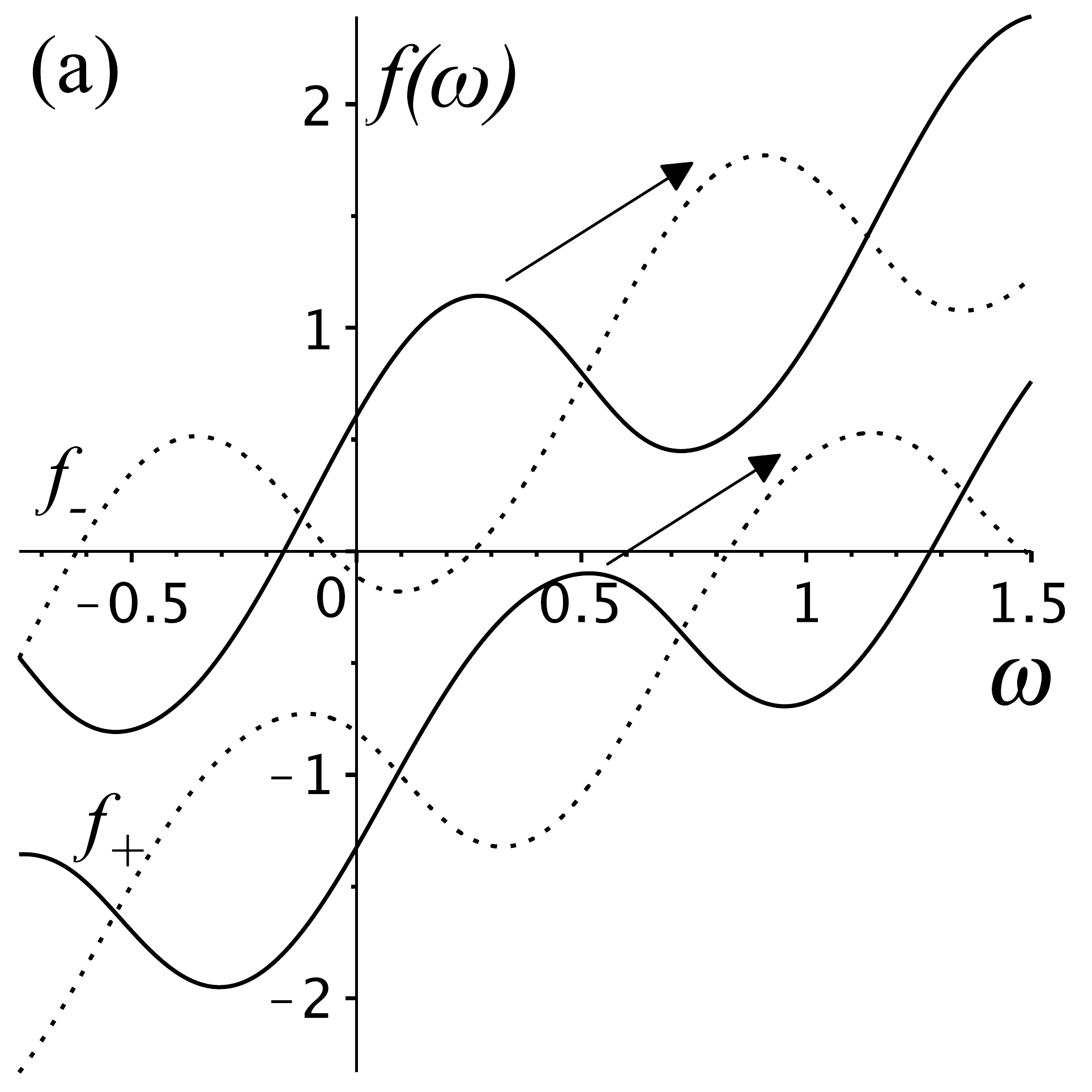}\hskip
1.5 cm\includegraphics[width=0.3\linewidth]{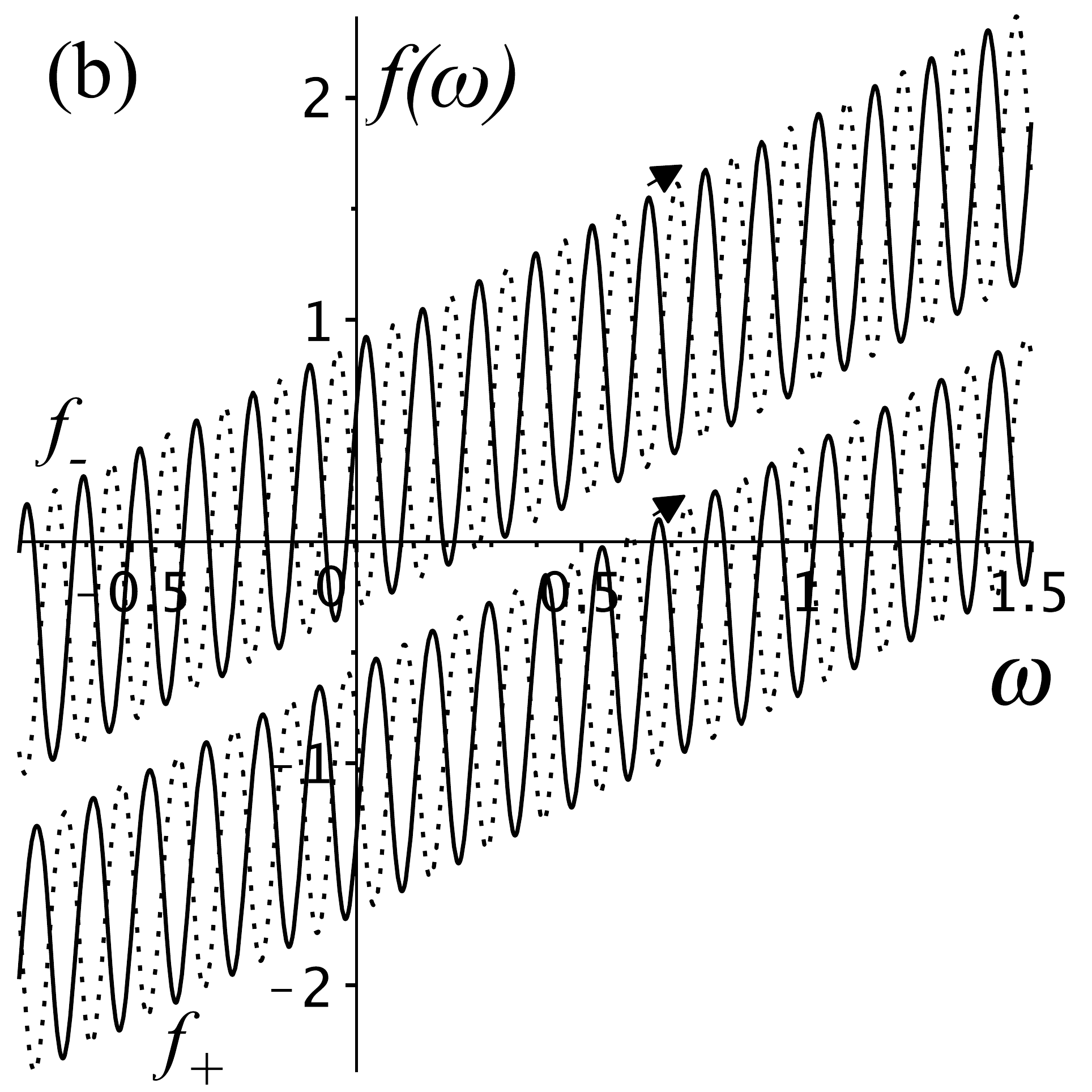}

\caption{Zeros of functions $f_{\pm}\left(\omega\right)$ determine frequencies
$\omega$ of plane waves. $f_{\pm}\left(\omega\right)$ are calculated
from Eq. \eqref{eq:omega} for (a) $\tau=5$ and (b) $\tau=50$. Solid
lines: $\varphi=0,$ dashed lines: $\varphi=\pi.$ Other parameters:
$q=0,$ $\delta=0.4,$ $\eta=0.5,$ $\beta=0.5,$ $\epsilon=1,$ $\mu=-1$,
and $\nu=-0.1$.\label{fig:pwd-sol}}
\end{figure}

\begin{figure}[h]
\centering \includegraphics[width=0.3\linewidth]{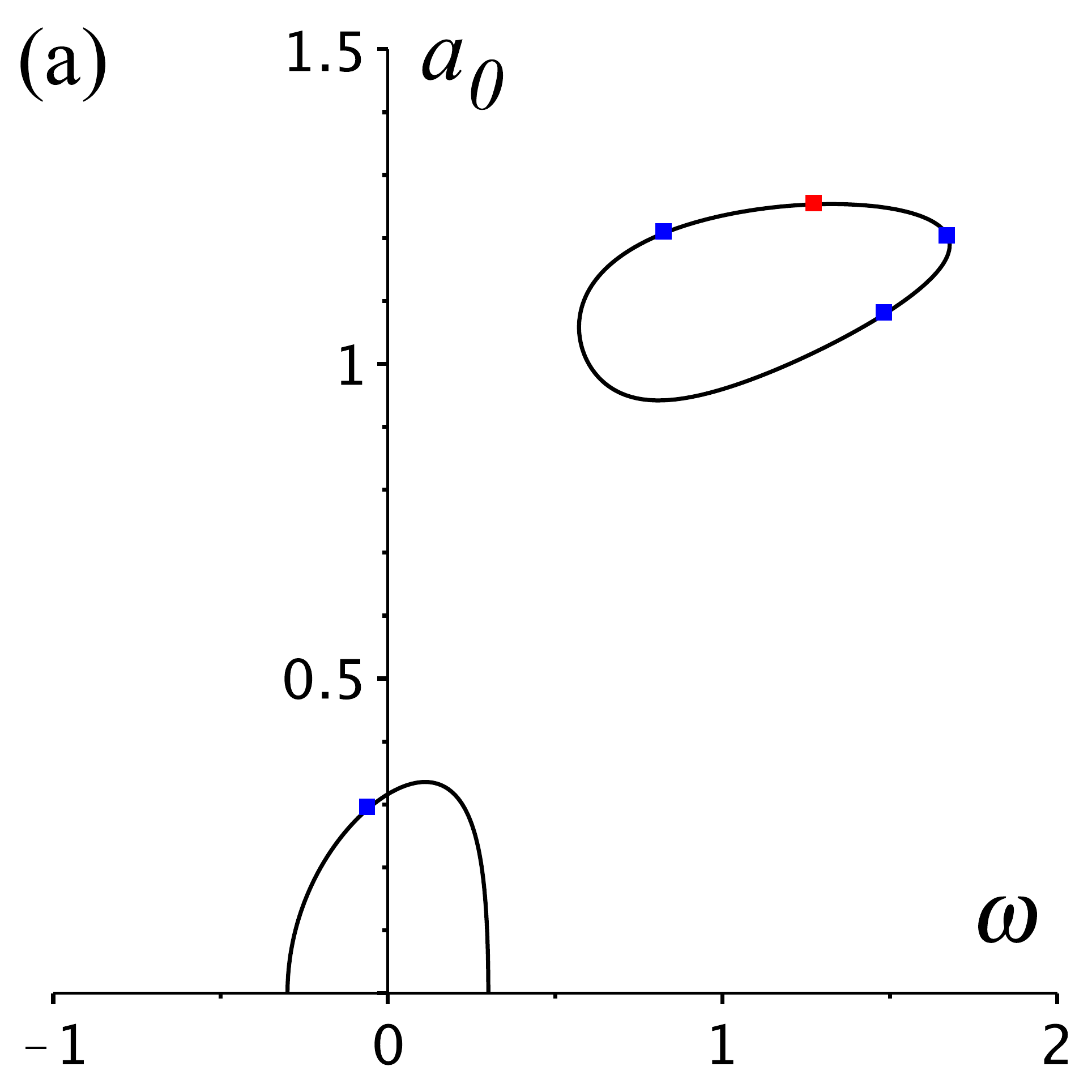}\hskip
1.5 cm\includegraphics[width=0.3\linewidth]{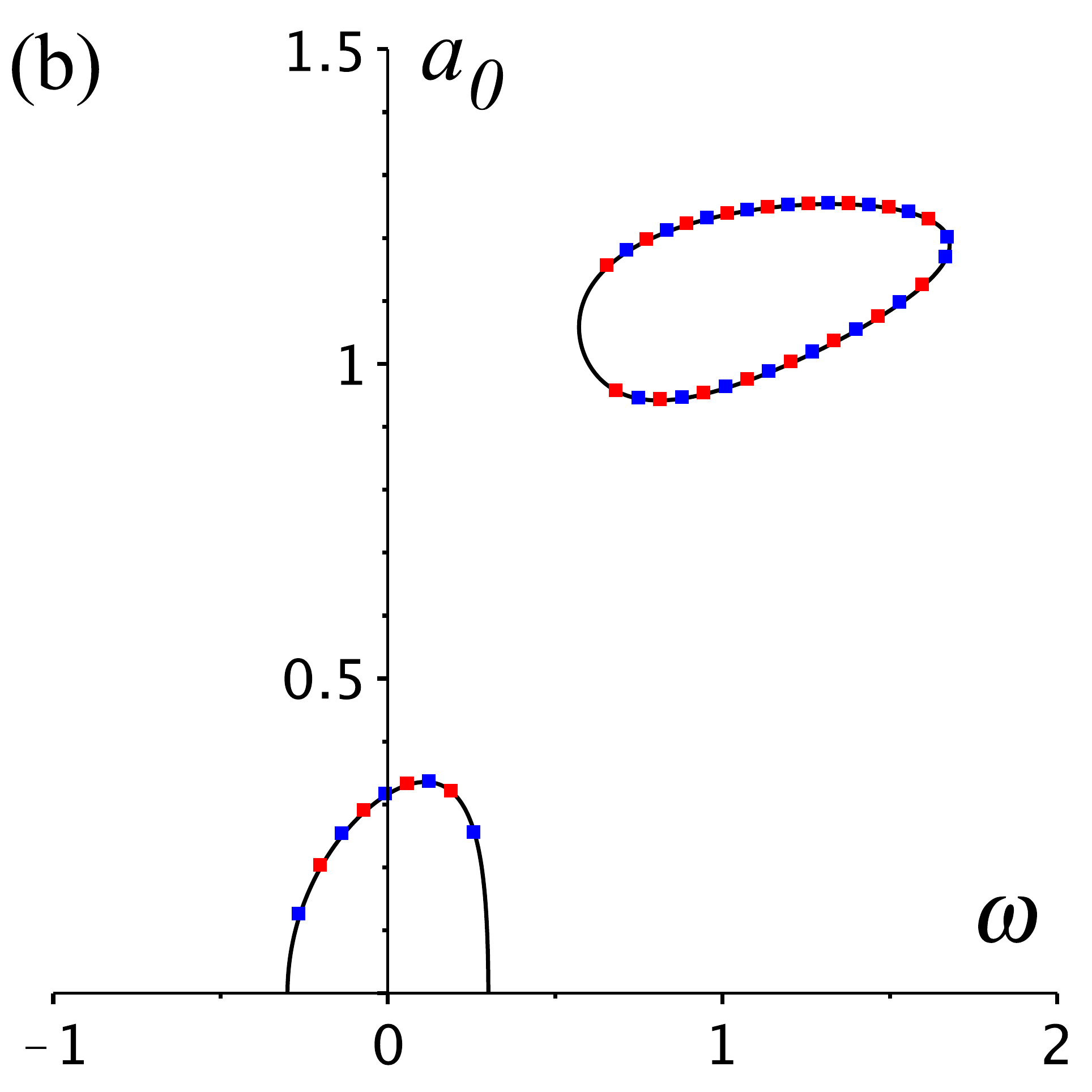}

\caption{Amplitudes $a_{0}$ and frequencies $\omega$ of plane wave solutions
for fixed $\delta=0.4$. The solutions are located on curves determined
by Eq. \eqref{eq:pwd-ell1} (solid lines), which are independent of
the delay $\tau$. Exact plane wave solutions for different values
of $\tau$ are shown by points: (a) $\tau=5$, (b) $\tau=50$. Red
points: $\varphi=0,$ blue points: $\varphi=\pi.$ Other parameters:
$q=0,$ $\eta=0.5,$ $\beta=0.5,$ $\epsilon=1,$ $\mu=-1$, and $\nu=-0.1$.\label{fig:pwd-ell2}}
\end{figure}

Solving \eqref{eq:pw-delay1} with respect to $a_{0}^{2}$ and substituting
the result into \eqref{eq:pw-delay2}, we obtain the equation for
the frequencies $\omega$ of the plane waves 
\begin{eqnarray}
0=f_{\pm}\left(\omega\right)=\omega+\frac{q^{2}}{2}+\frac{\epsilon\pm\sqrt{\epsilon^{2}-4\mu\left(\delta-\beta q^{2}+\eta\cos\left(\omega\tau-\varphi\right)\right)}}{2\mu}-\label{eq:omega} \\
-\frac{\nu\left(\epsilon\pm\sqrt{\epsilon^{2}-4\mu\left(\delta-\beta q^{2}+\eta\cos\left(\omega\tau-\varphi\right)\right)}\right)^{2}}{4\mu^{2}}+\eta\sin\left(\omega\tau-\varphi\right),\nonumber
\end{eqnarray}
which can be studied numerically for any fixed value of the delay
$\tau$. The functions $f_{\pm}(\omega)$ are shown on Fig.~\ref{fig:pwd-sol}.
In particular, with the increase of $\tau$ the functions $f_{\pm}\left(\omega\right)$
oscillate faster and number of solutions grows. This corresponds to
a general result obtained in \cite{PhysRevE.79.046221} stating that
in the limit of large delay the number of periodic solutions increases
linearly with $\tau$. Eventually the solutions fill the curves defined
by relation \eqref{eq:pwd-ell1}. Figure \eqref{fig:pwd-ell2} shows
the exact solutions $(a_{0},\omega)$ for increasing delay time $\tau$
as points on these curves. It is interesting to remark the strong
analogy between the curve \eqref{eq:pwd-ell1} of the plane wave solutions
and the ellipse of external cavity modes appearing in rate equation
models for semiconductor lasers with time delayed feedback, where
the notion is frequently used \cite{Mork1992,Heil2001,Yanchuk2010a,Soriano2013,Rottschafer2007}.

Let us consider the effect of the feedback phase $\varphi$, which
shifts the function $f_{\pm}\left(\omega\right)$, see Fig.~\ref{fig:pwd-sol}.
For small delay times, when $f_{\pm}\left(\omega\right)$ oscillates
slowly, the role of $\varphi$ is pronounced since a few of individual
solutions are moved along the curve $a_{0}(\omega)$ and new solutions
can be born. As $\tau$ increases, the overall number of plane wave
solutions becomes large, thus diminishing the effect of the feedback
phase, as shown in Figs.~\ref{fig:pwd-sol}(b) and \ref{fig:pwd-ell2}(b).

\begin{figure}[h]
\centering \includegraphics[width=0.3\linewidth]{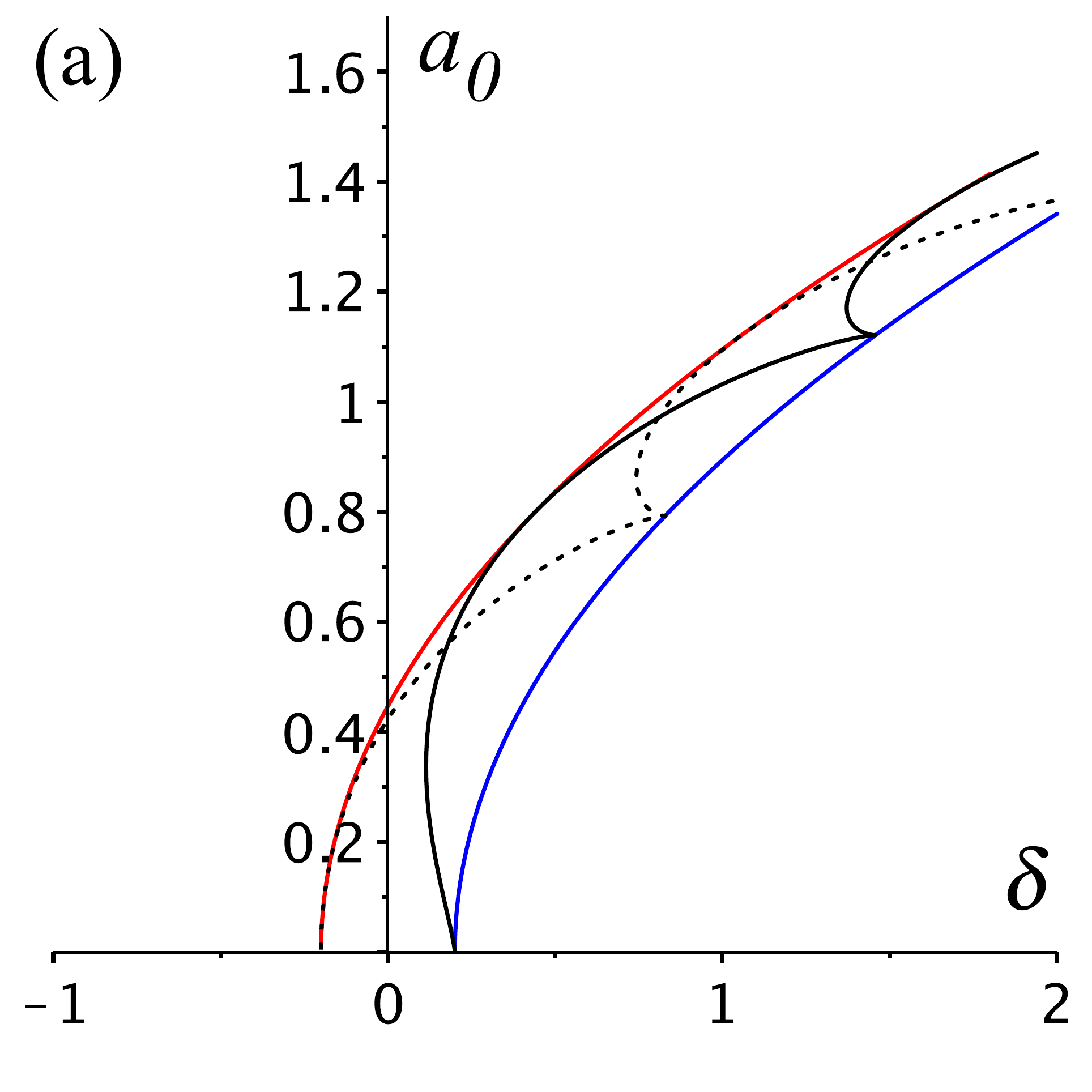}\hskip
1.5 cm\includegraphics[width=0.3\linewidth]{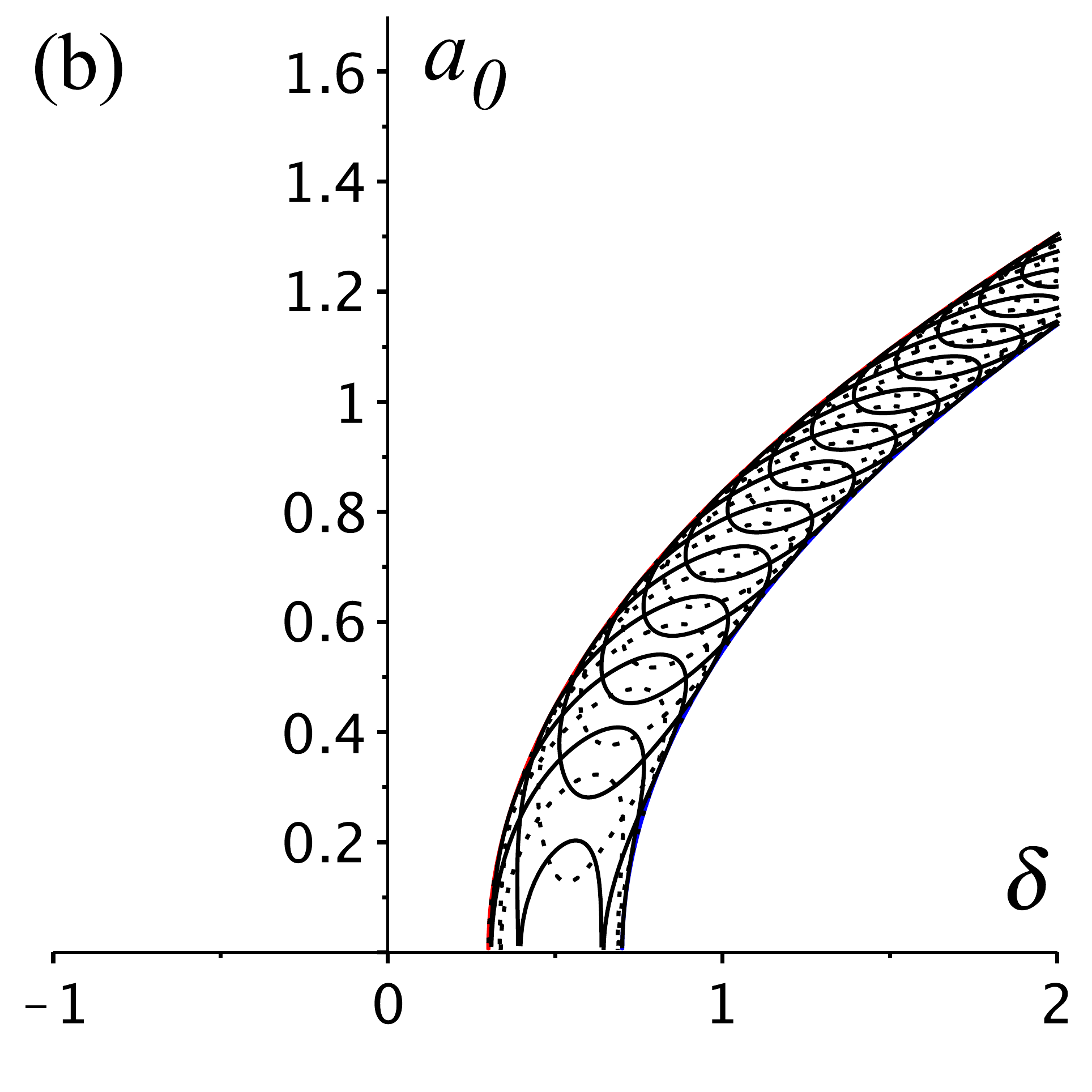}

\centering \includegraphics[width=0.3\linewidth]{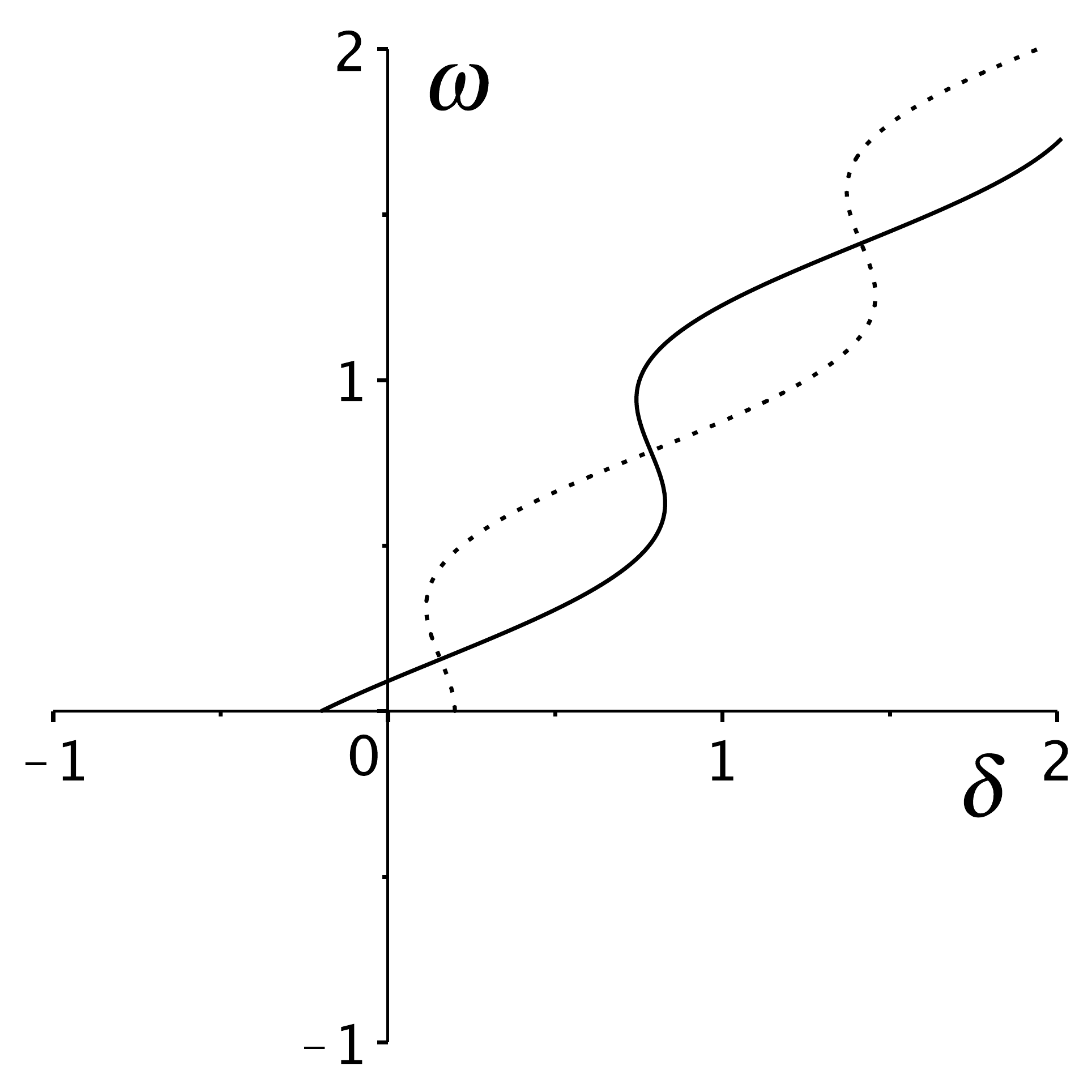}\hskip
1.5 cm\includegraphics[width=0.3\linewidth]{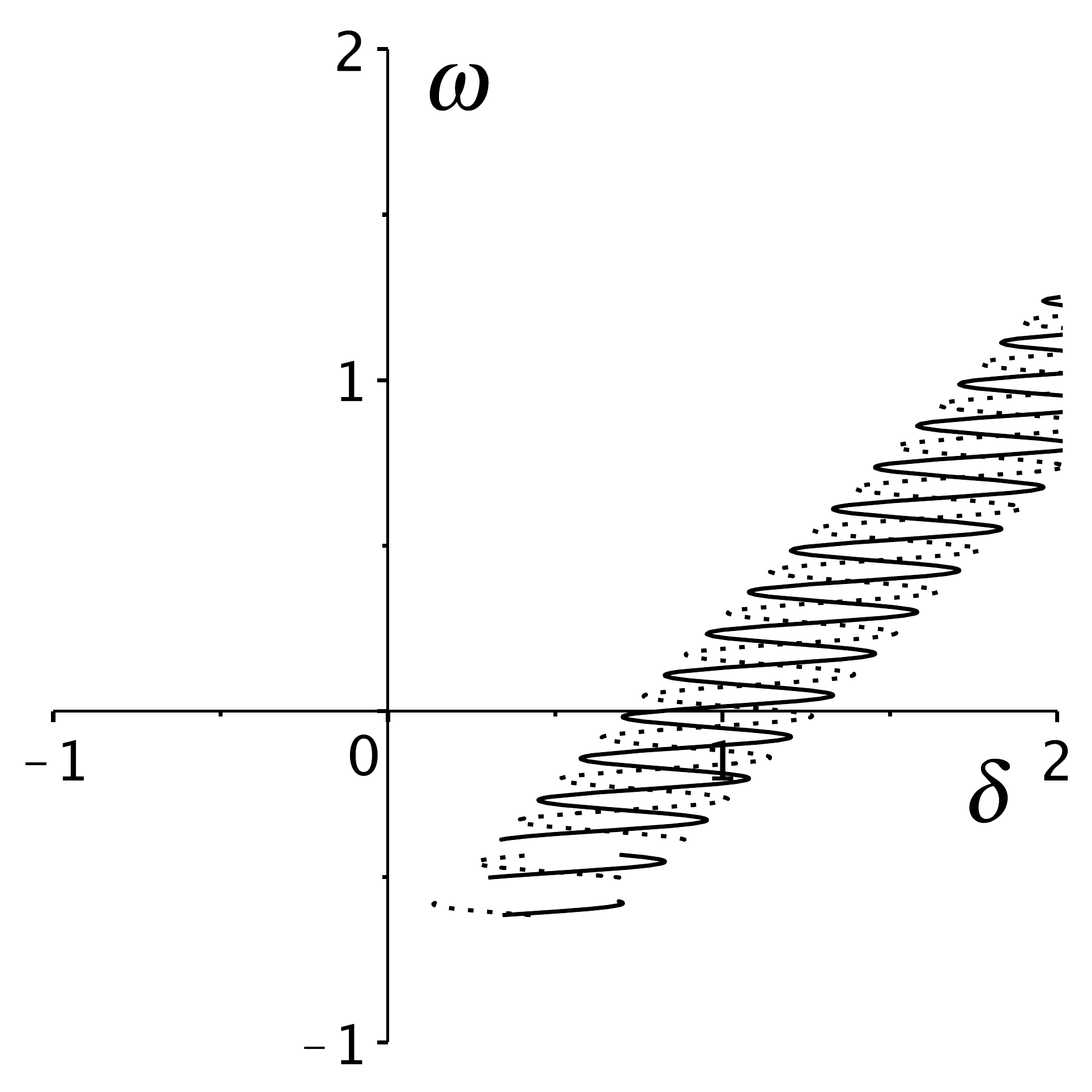}

\caption{Amplitudes $a_{0}\left(\delta\right)$ and frequencies $\omega\left(\delta\right)$
of plane wave solutions in cubic CGLE with delayed feedback for (a)
$\tau=5$, $q=0,$ (b) $\tau=50$, $q=1$. Solid black lines: $\varphi=0,$
dotted lines: $\varphi=\pi.$ The enveloping red and blue lines are
the branches of plane waves solutions for $\tau=0$, $\varphi=\pi$
and $\tau=0$, $\varphi=0$, respectively. Other parameters: $\eta=0.2,$
$\beta=0.5,$ $\epsilon=1$, $\mu=\nu=0$.\label{fig:pwd-sol_cubic}}
\end{figure}

\begin{figure}[h]
\centering \includegraphics[width=0.3\linewidth]{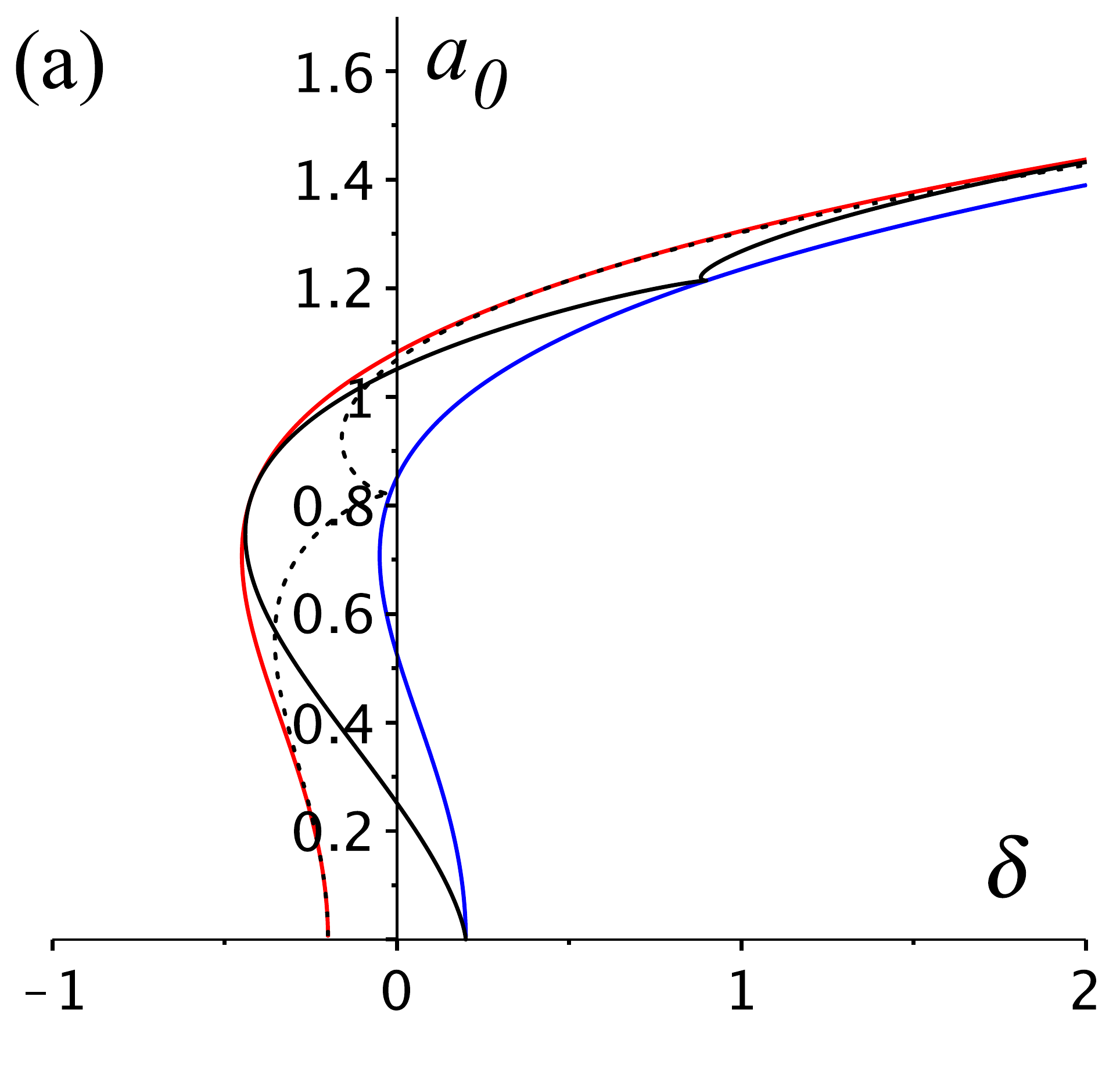}\hskip
1.5 cm\includegraphics[width=0.3\linewidth]{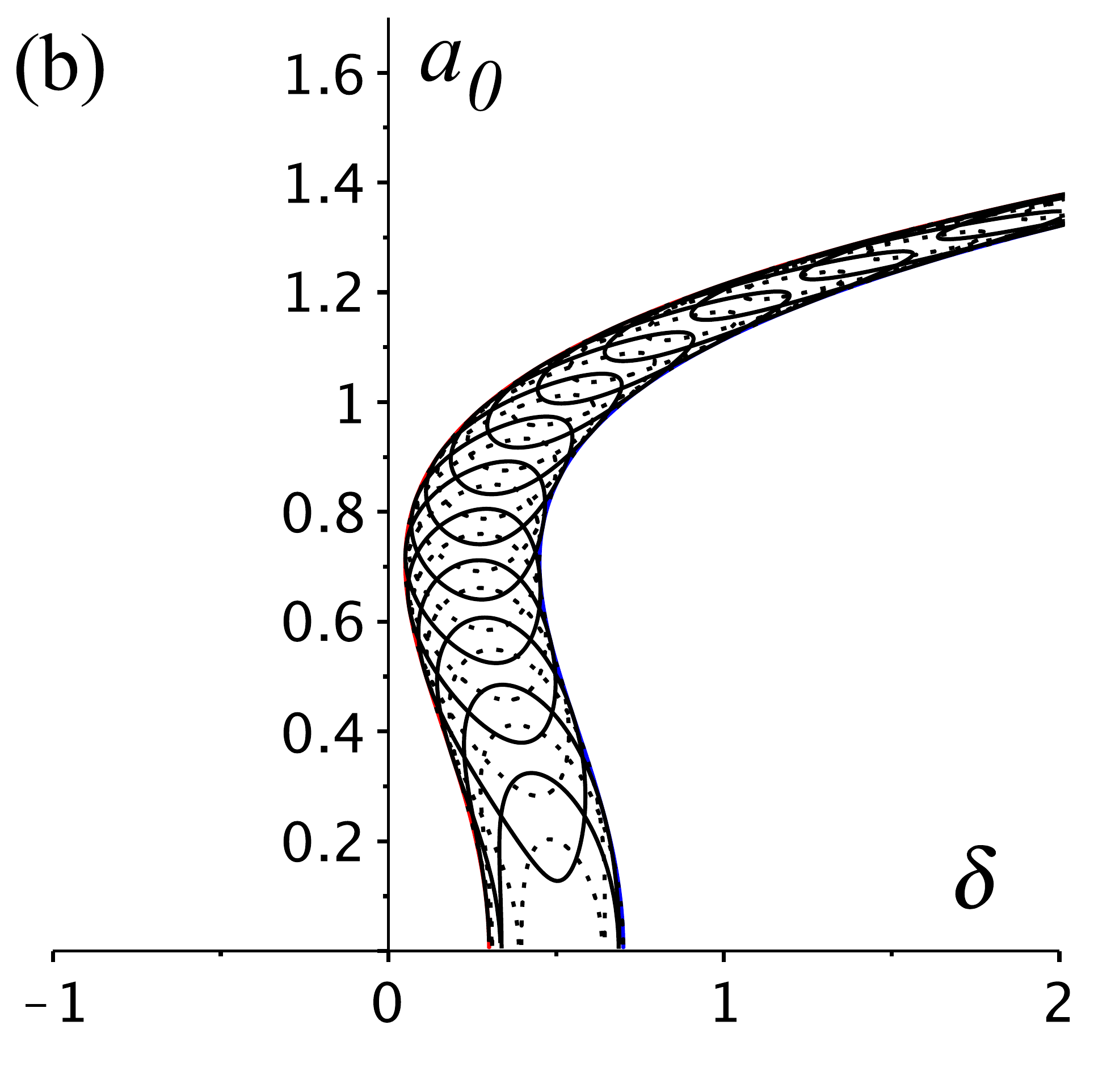}

\centering \includegraphics[width=0.3\linewidth]{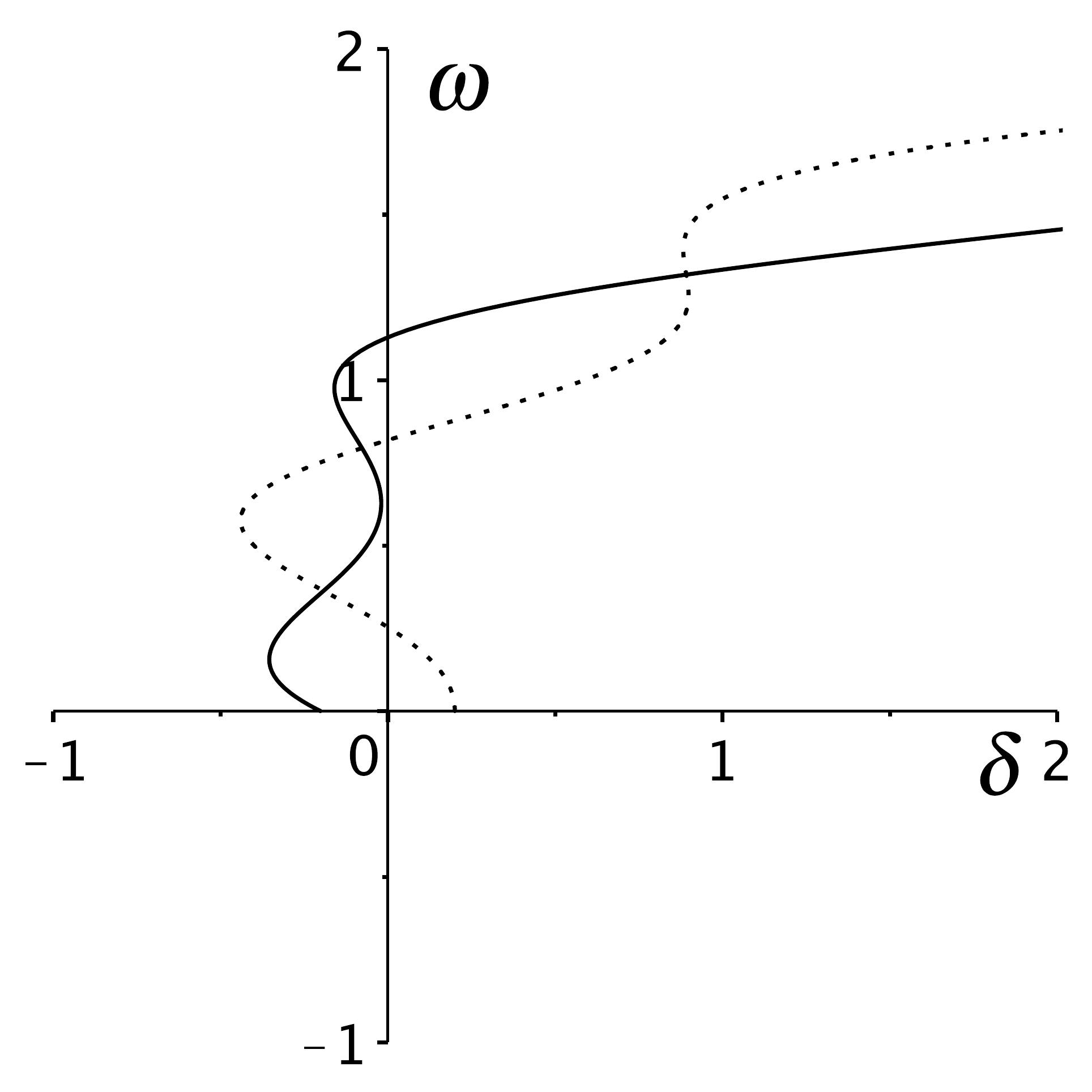}\hskip
1.5 cm\includegraphics[width=0.3\linewidth]{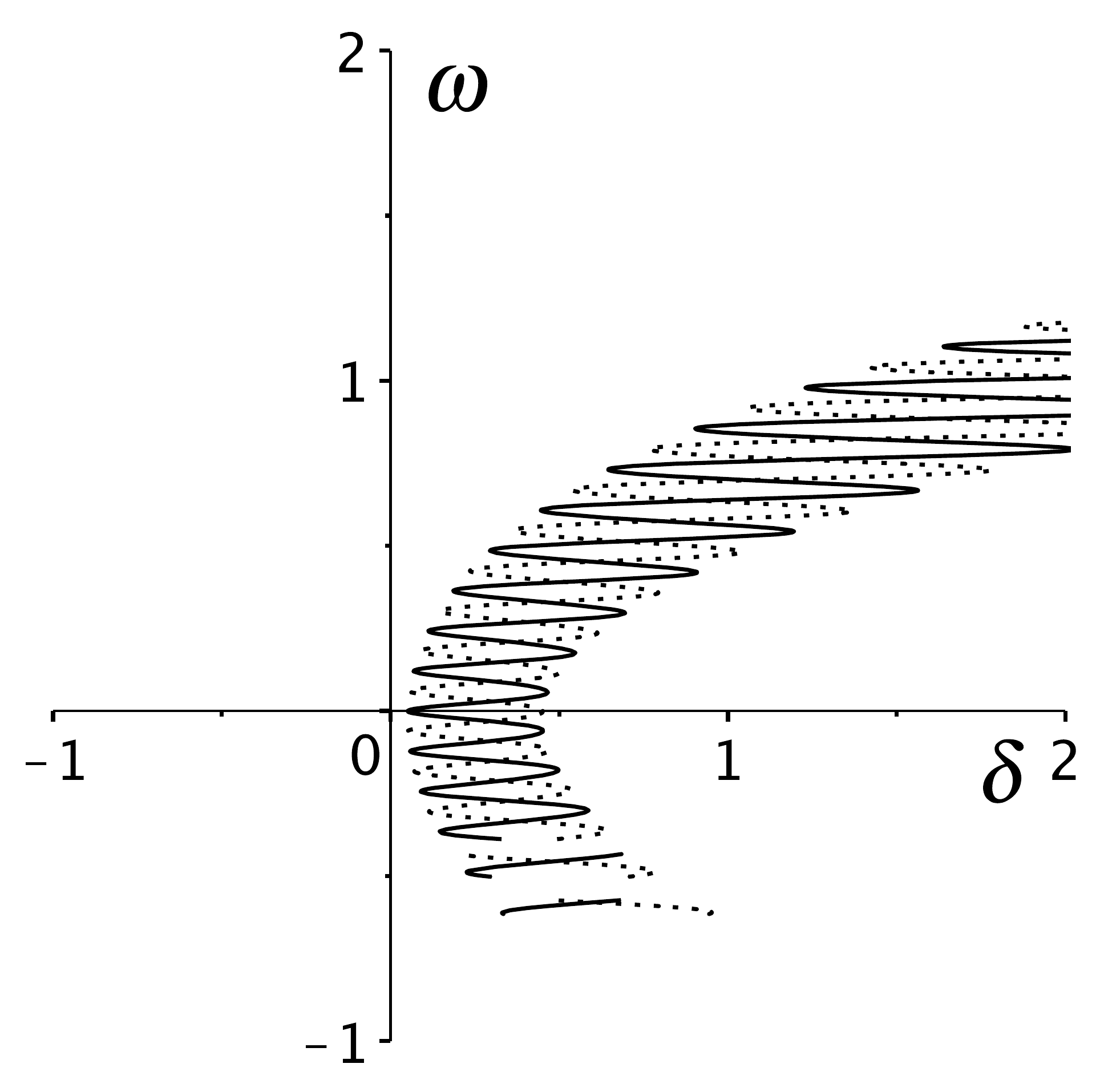}

\caption{Amplitudes $a_{0}\left(\delta\right)$ and frequencies $\omega\left(\delta\right)$
of plane wave solutions in quintic CGLE with delayed feedback for
(a) $\tau=5$, $q=0,$ (b) $\tau=50$, $q=1$. Solid black lines:
$\varphi=0,$ dotted lines: $\varphi=\pi.$ The enveloping red and
blue lines are the branches of plane wave solutions for $\tau=0$,
$\varphi=\pi$ and $\tau=0$, $\varphi=0$, respectively. Other parameters:
$\eta=0.2,$ $\beta=0.5,$ $\epsilon=1,$ $\mu=-1,$ $\nu=-0.1$.\label{fig:pwd-sol_quintic}}
\end{figure}

Even though exact values of $a_{0}$ and $\omega$ for any fixed set
of parameters can only be computed numerically, the branches of plane
wave solutions versus parameter $\delta$ can be obtained explicitly
in a parametric form $\left(a_{0}\left(\omega\right),\delta\left(\omega\right)\right)$.
Namely, the amplitude $a_{0}(\omega)$ is the solution of Eq. \eqref{eq:pw-delay2}.
Substituting $a_{0}(\omega)$ into Eq. \eqref{eq:pw-delay1} and solving
the resulting equation for $\delta$, we obtain an expression for
$\delta(\omega)$. As a result, the branches of plane waves have the
parametric form 
\begin{eqnarray}
a_{0}(\omega): & =\left\{ \mbox{solution of }\,\nu\left(a_{0}^{2}\right)^{2}+a_{0}^{2}-\eta\sin(\omega\tau-\varphi)-\omega-\frac{q^{2}}{2}=0\right\} ,\label{eq:branch}\\
\delta(\omega): & =\beta q^{2}-\epsilon a_{0}^{2}(\omega)-\mu a_{0}^{4}(\omega)-\eta\cos(\omega\tau-\varphi).\nonumber 
\end{eqnarray}
The branches of plane wave solutions obtained using this procedure
are shown in Fig.~\ref{fig:pwd-sol_cubic} and Fig.~\ref{fig:pwd-sol_quintic}
for cubic and quintic CGLE, respectively. Interestingly, the branches
have the form of snaking curves, which are constrained between two
``limiting'' branches \eqref{eq:wsol1} (red and blue curves in
Figs.~\ref{fig:pwd-sol_cubic}-\ref{fig:pwd-sol_quintic}), which
can be obtained by setting $\tau=0$, $\varphi=0$ and $\tau=0$,
$\varphi=\pi$, respectively. One can see that the increase of $\tau$
leads to even more dense snaking of the curve (Fig. \ref{fig:pwd-sol_cubic}(b)),
whereas the increase of $q$ shifts the branches to the right. 

\subsubsection*{Stability of plane wave solutions}

As in the case without delayed feedback, we use the ansatz \eqref{eq:pw-ansatz}
to investigate the stability of plane waves. After substituting it
into delayed CGLE \eqref{CGLE} and linearization in $a_{p}$, we
obtain: 
\begin{eqnarray}
\partial_{t}a_{p}+i\omega\left(a_{0}+a_{p}\right)=\left(\beta+\frac{i}{2}\right)\left(\partial_{xx}a_{p}+2iq\partial_{x}a_{p}-q^{2}(a_{0}+a_{p})\right)+\delta\left(a_{0}+a_{p}\right)+\nonumber \\
+\left(\epsilon+i\right)\left(a_{0}^{3}+2a_{0}^{2}a_{p}+a_{0}^{2}\overline{a_{p}}\right)+\left(\mu+i\nu\right)\left(a_{0}^{5}+3a_{0}^{2}a_{p}+2a_{0}^{2}\overline{a_{p}}\right)\label{eq:pw2-1}\\
+\eta e^{i\varphi-i\omega\tau}\left(a_{0}+a_{p}(t-\tau)\right).\nonumber 
\end{eqnarray}
We simplify this equation using Eq. \eqref{eq:pwdelaymain} and substitute
the exponential ansatz \eqref{eq:ap} for $a_{p}(x,t)$ into it. As
a result, similarly to Section \ref{sub:Stability-of-pw}, we obtain
a linear system of equations with respect to two unknowns $a_{+}$
and $a_{-}$: 
\[
M_{\tau}\left(\begin{array}{c}
a_{+}\\
a_{-}
\end{array}\right)=0,
\]
where{\scriptsize 
\begin{equation}
M_{\tau}=\begin{bmatrix}\lambda+i\omega-\delta+\left(\beta+\frac{i}{2}\right)\left(k^{2}+2kq+q^{2}\right)- & \qquad & -\left(\epsilon+i\right)a_{0}^{2}-2\left(\mu+i\nu\right)a_{0}^{4}\\
-2\left(\epsilon+i\right)a_{0}^{2}-3\left(\mu+i\nu\right)a_{0}^{4}-\eta e^{-\lambda\tau}e^{i\varphi-i\omega\tau}\\
\\
 &  & \lambda-i\omega-\delta+\left(\beta-\frac{i}{2}\right)\left(k^{2}-2kq+q^{2}\right)-\\
-\left(\epsilon-i\right)a_{0}^{2}-2\left(\mu-i\nu\right)a_{0}^{4} &  & -2\left(\epsilon-i\right)a_{0}^{2}-3\left(\mu-i\nu\right)a_{0}^{4}-\eta e^{-\lambda\tau}e^{-i\varphi+i\omega\tau}
\end{bmatrix}.
\end{equation}
}The condition 
\begin{equation}
\det M_{\tau}=0\label{eq:detM}
\end{equation}
now gives us the characteristic equation for the perturbation growth
rate $\lambda$.

\begin{figure}[h]
\centering \includegraphics[width=0.3\linewidth]{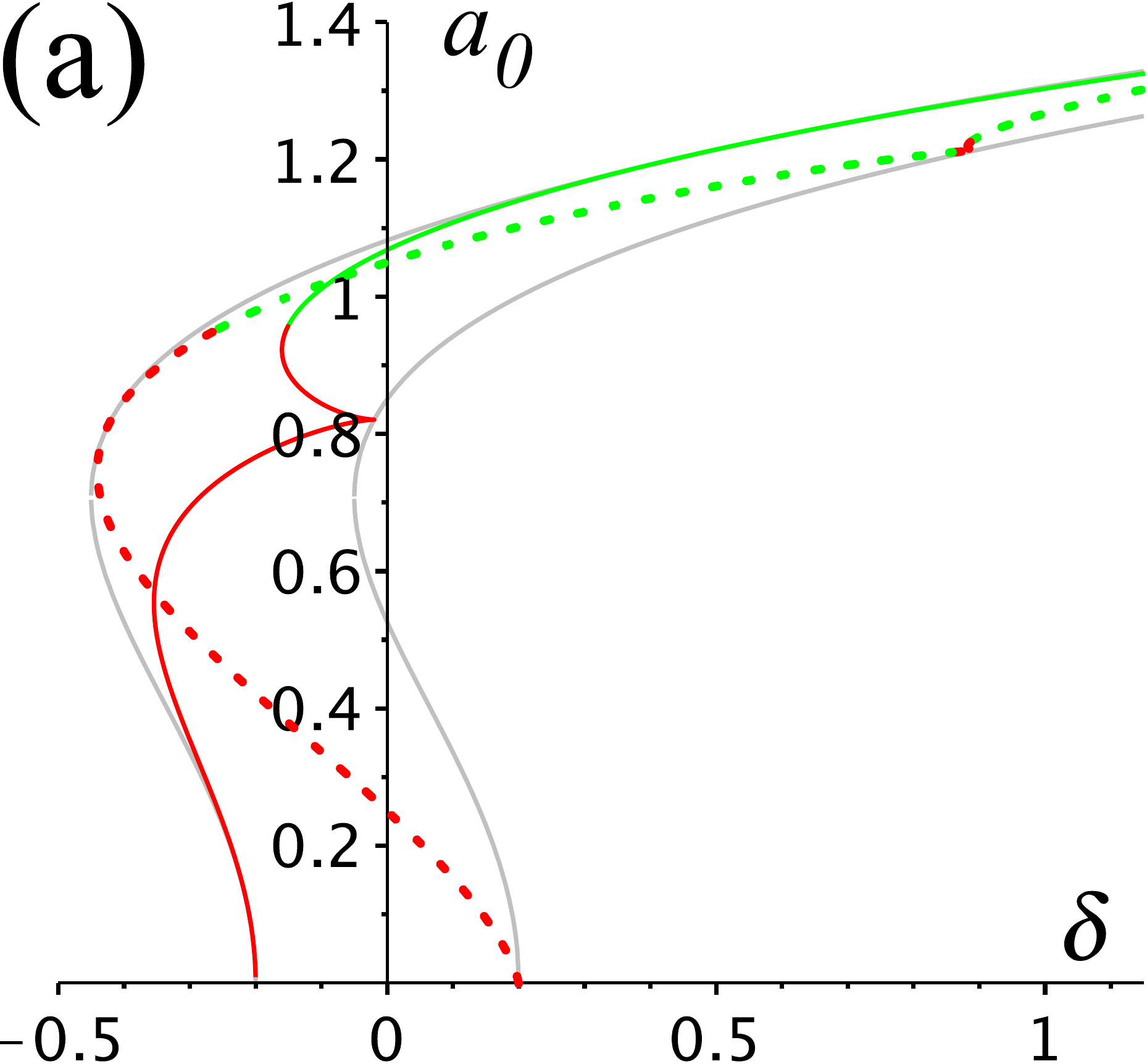}\hskip
1.5 cm\includegraphics[width=0.3\linewidth]{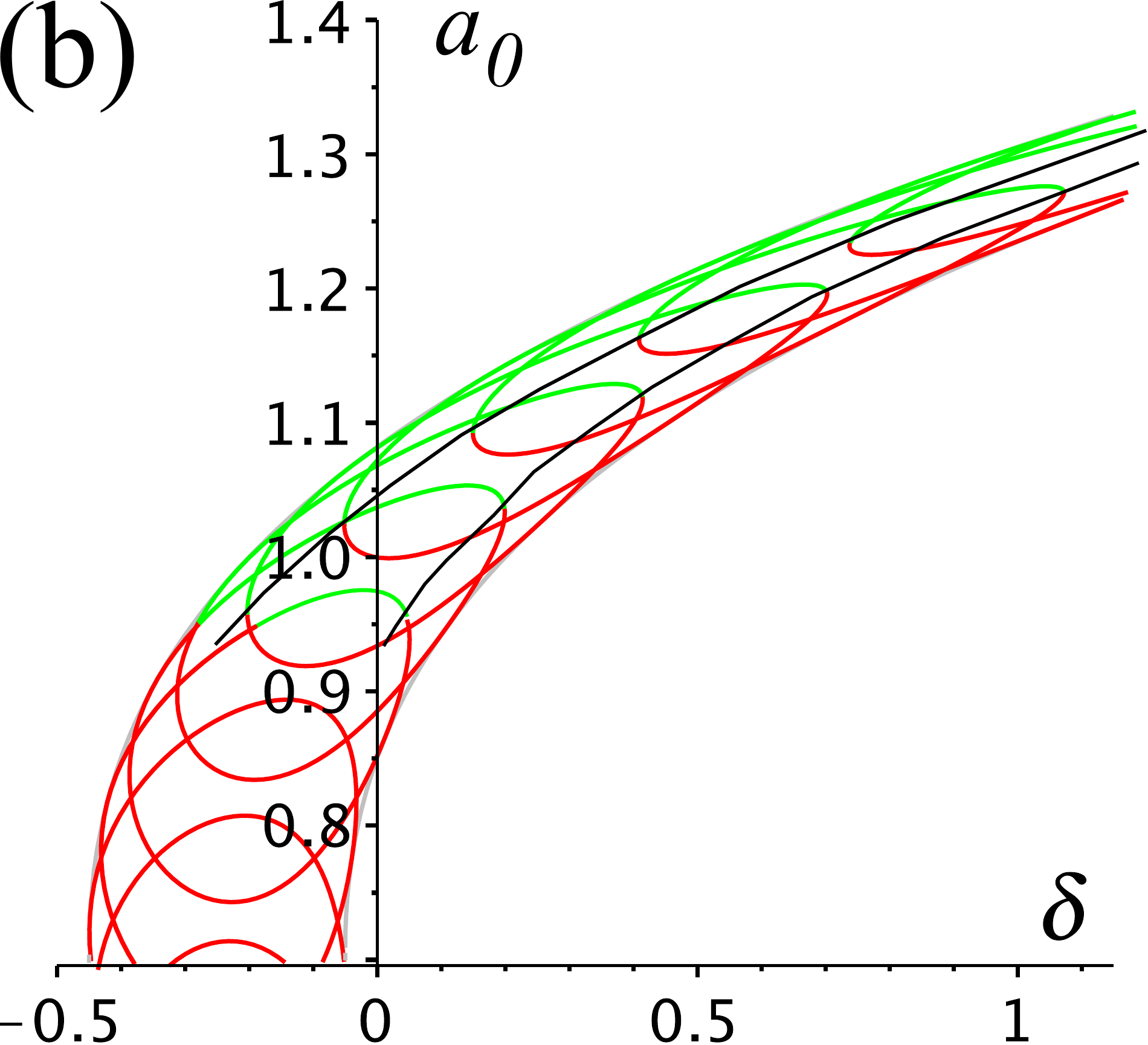}

\caption{Stability of the branches of plain wave solutions for different delay
times (a) $\tau=5,\:\varphi=0$ (solid line), $\varphi=\pi$ (dotted
line) (b) $\tau=50,\:\varphi=0$. Stable parts are shown in green
and unstable -- in red. Black lines in (b) denote the destabilization
borders in the limit of large delay. Other parameters: $\eta=0.2,$
$\beta=0.5,$ $\epsilon=1,$ $\mu=-1$, and $\nu=-0.1$.\label{fig:pwd-stab-q0}}
\end{figure}

Stability of individual plane wave solutions with arbitrary delay
time $\tau$ is determined by the real parts of the roots of the characteristic
equation \eqref{eq:detM}. If for all $k$, the roots satisfy $\Re[\lambda]<0$
(except the trivial one with $\Re[\lambda]=0$), then the plane wave
is stable. For a fixed $q$, the branch of plane waves is given parametrically
by $\left(a_{0}\left(\omega\right),\delta\left(\omega\right)\right)$,
defined by Eq. \eqref{eq:branch}. Substituting $a_{0}(\omega)$ and
$\delta(\omega)$ in Eq. \eqref{eq:detM}, we obtain a nonlinear characteristic
equation for $\lambda$ 
\begin{equation}
{\cal F}\left(\lambda;q,\omega,k,\mathbf{p}\right)=0,\label{eq:F}
\end{equation}
where $\mathbf{p}$ denotes system parameters $(\beta,\epsilon,\mu,\nu,\eta,\varphi,\tau)$.
Equation \eqref{eq:F} was solved numerically for fixed values of
$\mathbf{p}$ and varying $\omega$. In this way we obtain the stability
properties for the parametrically defined family of plane waves $\left(a_{0}\left(\omega\right),\delta\left(\omega\right)\right)$
for the given CGLE parameters and wavenumber $q.$ Figure \ref{fig:pwd-stab-q0}
shows the stability properties on the branches of plane waves $\left(a_{0}\left(\omega\right),\delta\left(\omega\right)\right)$
with $q=0$ and various delay times $\tau$. Stable solutions are
plotted in green, unstable in red. Figure \ref{fig:pwd-stab-q0}(a)
illustrates the effect of the feedback phase $\varphi$ on plane wave
solutions and their stability at relatively small delay time $\tau=5$,
whereas Fig.~\ref{fig:pwd-stab-q0}(b) depicts the bifurcations for
larger delay time $\tau=50$. One can observe the growing number of
multistable plane waves with the increase of the delay. With the increase
of the parameter $\delta$ high amplitude parts of the snaking branches
with the higher amplitude becomes stable while low amplitude parts
remain unstable. An additional analytical insight into the structure
of stable and unstable regions is obtained by using the large delay
approximation, which is discussed in the next section.

\subsection{Large delay\label{sub:The-case-of-large}}

For large delay times, the plane wave solutions fill the tube defined
by Eq. \eqref{eq:pwd-ell1} densely, see also Figs.~\ref{fig:soltube-1}
and \ref{fig:pwd-ell2}(b). Hence, instead of looking at individual
solutions and solving the transcendental Eq. \eqref{eq:omega}, it
is convenient to parametrize the whole family of solutions by a parameter
$\theta=(\omega\tau\text{-\ensuremath{\varphi}+\ensuremath{\pi}})\,\mbox{mod}\,2\pi$,
which may be represented as an angular coordinate on the tube. Every
single solution for given control parameters and wavevector $q$ can
be uniquely defined by the coordinate $\theta$. Therefore, one can
consider the question about the stability of a plane wave at a given
$\theta$-value with the amplitude $a_{0}(\theta)$ and the frequency
$\omega(\theta)$. The growth rate $\lambda$, which determines the
stability of an individual plane wave solution, is obtained from the
characteristic equation depending just on the system parameters and
the coordinate $\theta.$ Note that $\theta=\pi$ corresponds to the
plane waves with the maximal amplitude $a_{0}$. The two sides of
the tube of plane wave solutions, one with $\theta<\pi$ and another
with $\theta>\pi$, are projected onto the same set in the $(\delta,a_{0})$
plane, but correspond to different values of $\omega.$

\begin{figure}[h]
\centering \includegraphics[width=0.3\linewidth]{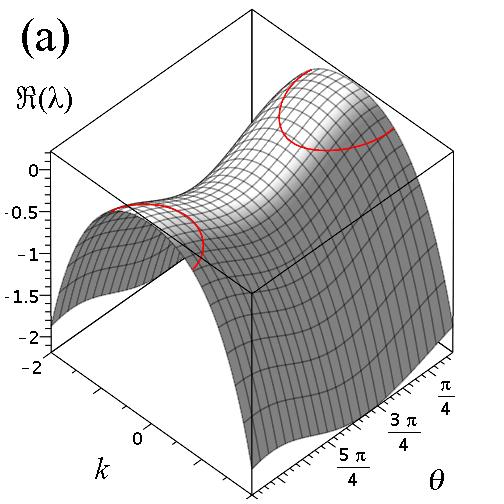}\hskip
1.5 cm\includegraphics[width=0.3\linewidth]{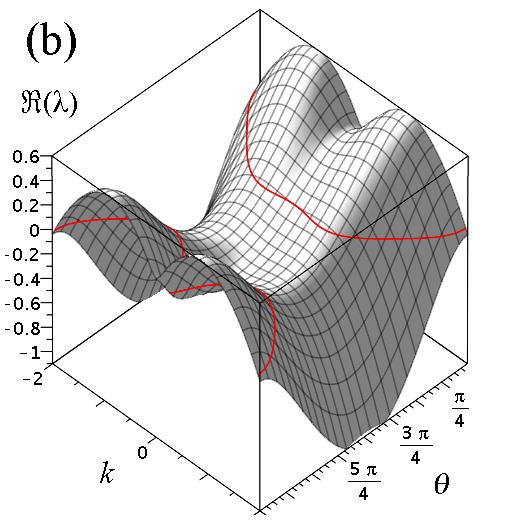}

\caption{Real parts of the eigenvalues $\Re[\lambda(k,\theta)]$ describing
the strong instability of plane waves. For a given plane wave, which
is determined by the parameter $\theta$, the positiveness of $\Re[\lambda(k,\theta)]$
for some $k$ implies the strong instability. Red curves depict zero
level lines. (a) $q=0,$ (b) $q=1$. On panel (a), the plane waves
are strongly unstable for $\theta\lesssim1.9$ and $\theta\gtrsim5.0$.
On panel (b), the plane waves are strongly unstable for $\theta\lesssim1.9$
and $\theta\gtrsim4.5$. Other parameters are: $\delta=0.5,\:\eta=0.2,$
$\beta=0.5,$ $\epsilon=1,$ $\mu=-1$, and $\nu=-0.1$.\label{fig:pwd-sol-str}}
\end{figure}

To study the strong instability of plane wave solutions, again, as
in Sec.~\ref{sub:The-case-of-large}, we neglect the terms containing
$e^{-\text{\ensuremath{\lambda}}\tau}$ in Eq. \eqref{eq:detM}. This
gives us a quadratic equation for $\lambda$.\textcolor{red}{{} }For
a given parameter $\theta$, the plane wave is strongly unstable,
if the maximum of $\Re[\lambda(k,\theta)]$ is positive. Figure \ref{fig:pwd-sol-str}
shows the the real part of $\lambda(k,\theta)$ for reduced equation
\eqref{eq:detM} for two different values of the wavenumber, $q=0$
and $q=1$ ($\delta$ is fixed to $0.5$). Red curves depict zero
level lines. We observe, that at larger values of $q$, unlike the
feedback-free case, the destabilization occurs first at non-zero values
of the perturbation wavenumber $k$.

\begin{figure}[h]
\centering \includegraphics[width=0.3\linewidth]{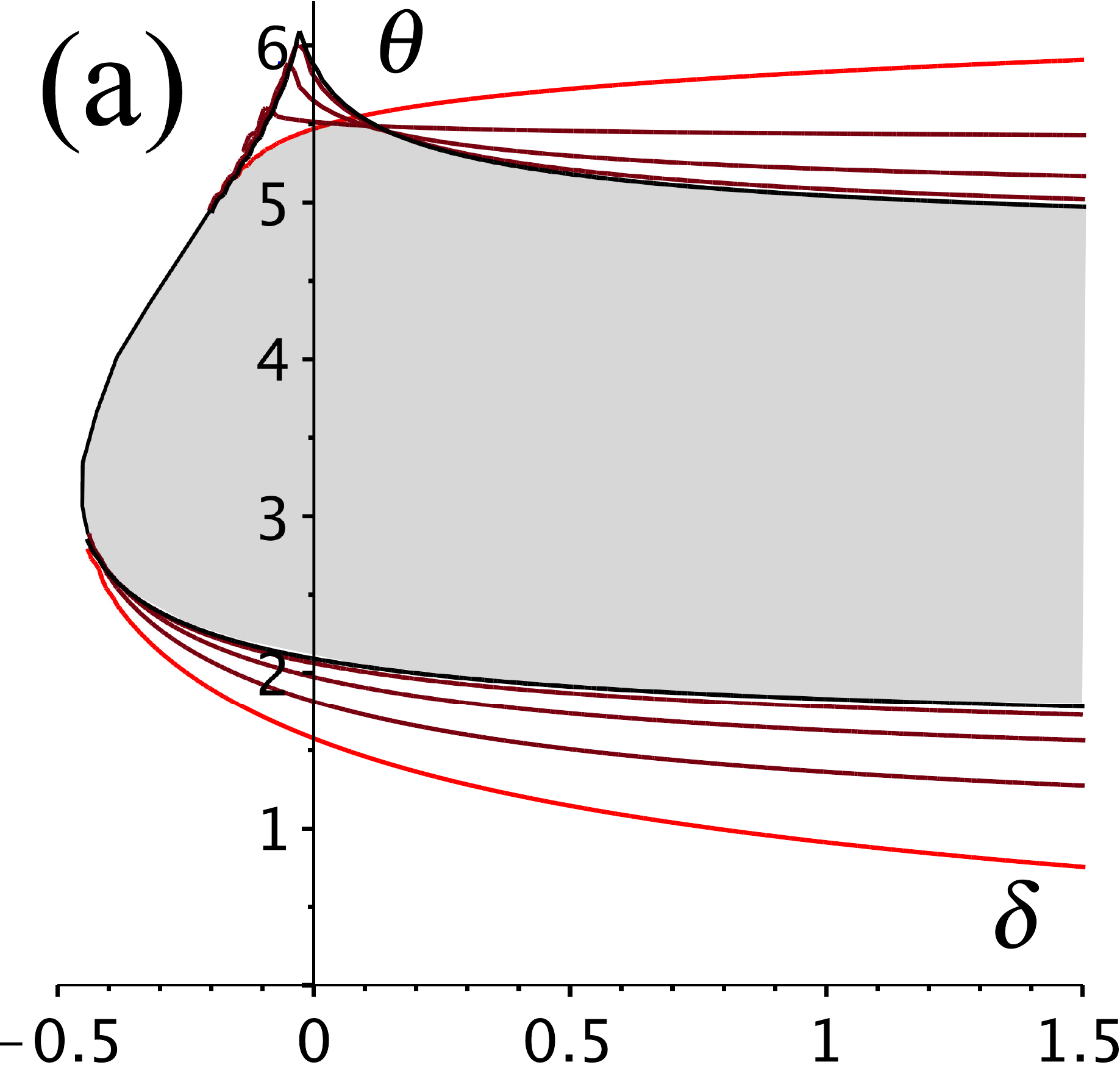}\hskip
1.5 cm\includegraphics[width=0.3\linewidth]{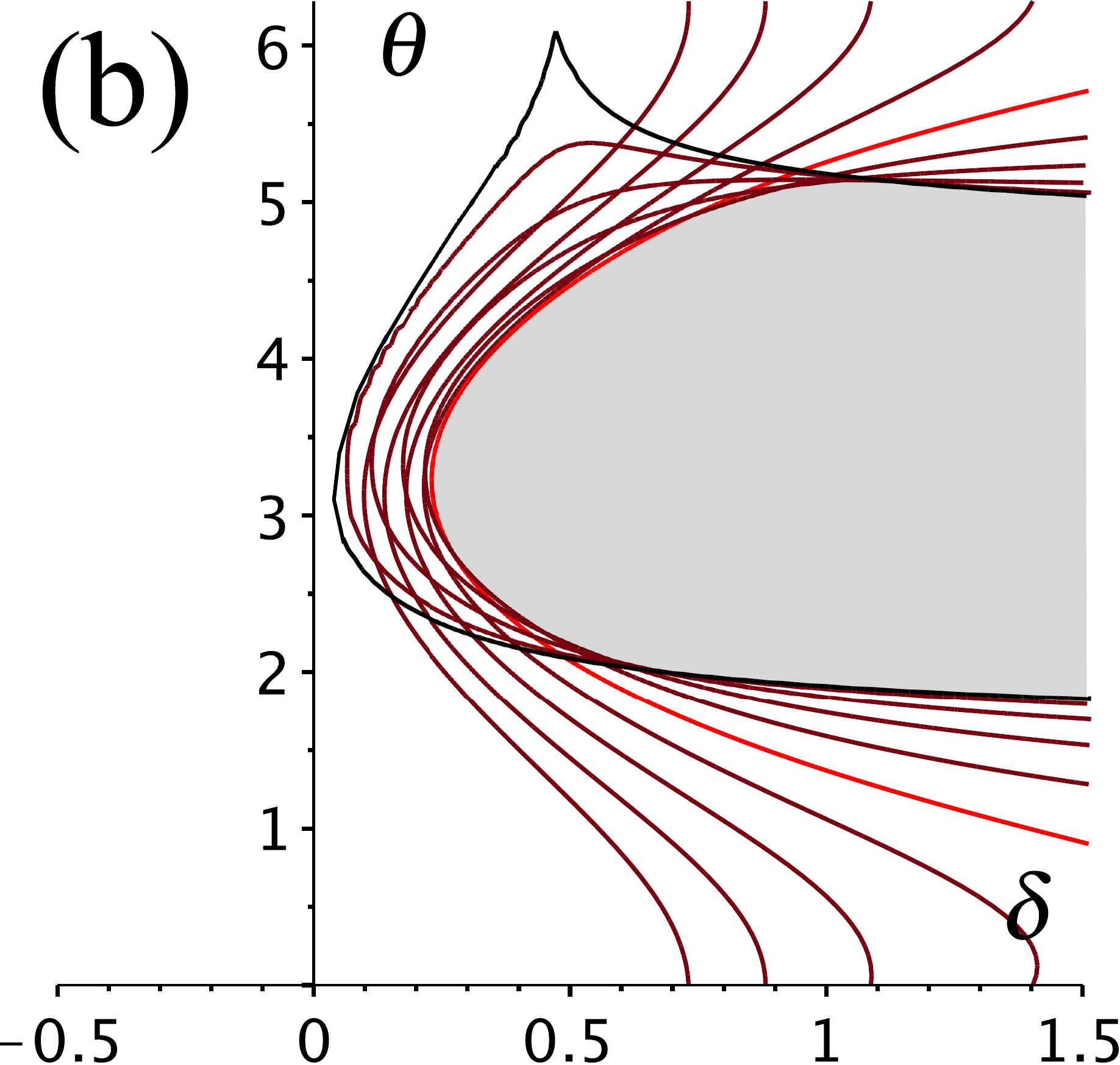}

\caption{Strong instability of the plane wave solutions with $q=0$ (a) and
$q=1$ (b).Gray region shows the absence of strong instability. Different
curves indicate zero contour levels of $\Re[\lambda(\delta,\theta)]$
corresponding to different perturbation wavenumbers $k.$ (a) $q=0,$
(b) $q=1$. Other parameters: $\eta=0.2,$ $\beta=0.5,$ $\epsilon=1,$
$\mu=-1$, and $\nu=-0.1$.\label{fig:pwd-bif-str}}
\end{figure}

In order to determine the strong instability boundary of the plane
wave solutions, in Fig. \ref{fig:pwd-bif-str} we plot on $(\delta,\theta)$
plane zero contour levels $\Re[\lambda(\delta,\theta)]=0$ corresponding
to different perturbation wavenumbers $k$. Figure \ref{fig:pwd-bif-str}(a)
shows that for the plane wave with $q=0$, stability border almost
coincides with the zero level line corresponding to $k\text{=}0$,
depicted by the black line. By contrast, for the plane wave with $q=1$
a significant part of stability border, shown in red color, corresponds
to destabilization with non-zero perturbation wavenumber $k\sim0.9$
(see Fig. \ref{fig:pwd-bif-str}(b)).

In the limit of large delay, the weak instability boundary is determined
by the pseudo-continuous spectrum 
\begin{equation}
\lambda=\frac{\gamma}{\tau}+i\xi.\label{PCS-1}
\end{equation}
To calculate this spectrum we substitute \eqref{PCS-1} into Eq. \eqref{eq:detM}
and neglect the terms proportional to $\gamma/\tau$. Then, denoting
\begin{equation}
Y=e^{-\gamma}e^{-i\xi\tau},\label{eq:PCS-2}
\end{equation}
one can solve the resulting quadratic equation with respect to $Y$.
The two roots $Y_{1}$ and $Y_{2}$ of this equation depend on the
model equation parameters $\mathbf{p}$, plane wave $(q,\theta)$,
wavenumber of the perturbation $k$, and the delay induced perturbation
modes $\xi$. For brevity we do not write down the explicit form of
these solutions. The factor $\gamma$ in \eqref{eq:PCS-2} is given
by 
\begin{equation}
\gamma_{1,2}=-\ln|Y_{1,2}(k,\xi;q,\theta,\mathbf{p})|.\label{eq:pw_weak_disp}
\end{equation}
For every plane wave solution defined by $q$ and $\theta$ and fixed
parameters $\mathbf{p}$ we obtain two surfaces $\gamma_{1,2}(k,\xi)$,
which generalize the dispersion relation (see Sec.~\eqref{sub:homlong}).\textcolor{red}{{}
}If $\gamma_{1,2}<0$ for all wavenumbers $k$ and delay modes $\xi$
(except the trivial eigenvalue corresponding to the Goldstone mode),
then the corresponding plane wave is stable (provided no strong instability
exists). Otherwise, it is weakly unstable. Thus, the weak instability
of the plane wave solutions is determined by considering the behavior
of the two surfaces $\gamma_{1,2}(k,\xi)$, or, more specifically,
the upper one.

\begin{figure}[h]
\centering \includegraphics[width=0.3\linewidth]{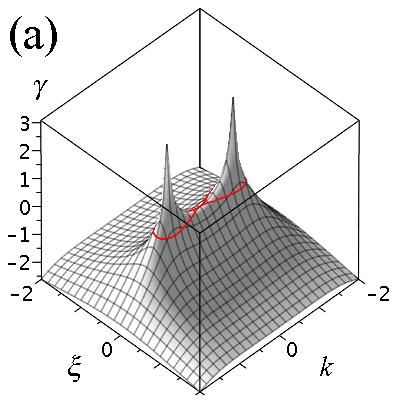}\hskip
0.5 cm\includegraphics[width=0.3\linewidth]{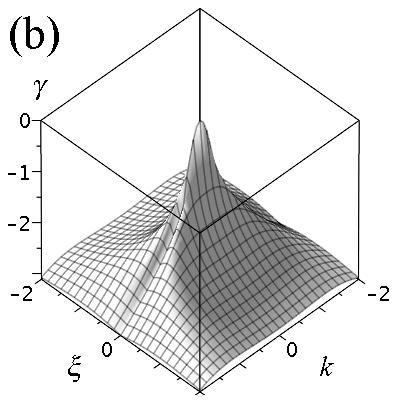}

\centering \includegraphics[width=0.3\linewidth]{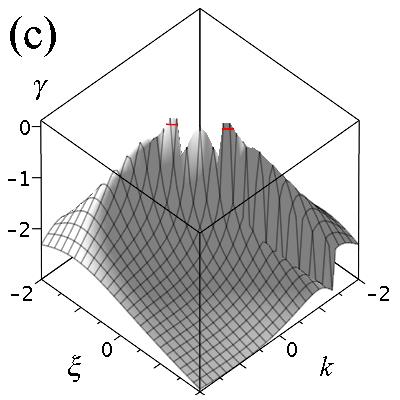}\hskip
0.5 cm\includegraphics[width=0.3\linewidth]{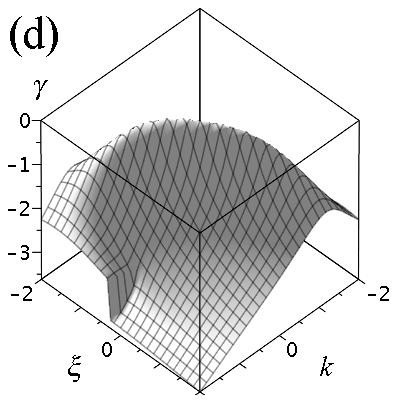}

\caption{The rescaled growth rate $\gamma=\tau\Re[\lambda]$ of different perturbation
modes (generalized dispersion relation) for different plane wave solutions.
$k$ is the spatial wavenumber of the perturbations, $\xi$ stands
for the delay induced modes. The positiveness of $\gamma$ implies
weak instability of the plane wave. The surface $\gamma(\xi,q)$ is
obtained from Eq. \eqref{eq:pw_weak_disp} for different values of
$\theta$, $\delta$, and $q$. (a) $\theta=0,\:\delta=0.5,\: q=0.$
(b) $\theta=2,\:\delta=1.0,\: q=0.$ (c) $\theta=2,\:\delta=1.0,\: q=1.$
(d) $\theta=4,\:\delta=1.0,\: q=1$. Panels (b) and (d) {[}(a) and
(c){]} correspond to stable {[}weakly unstable{]} plane waves. Other
parameters: $\eta=0.2,$ $\beta=0.5,$ $\epsilon=1,$ $\mu=-1$, and
$\nu=-0.1$.\label{fig:pwd-sol-wk}}
\end{figure}

\begin{figure}[h]
\centering \includegraphics[width=0.3\linewidth]{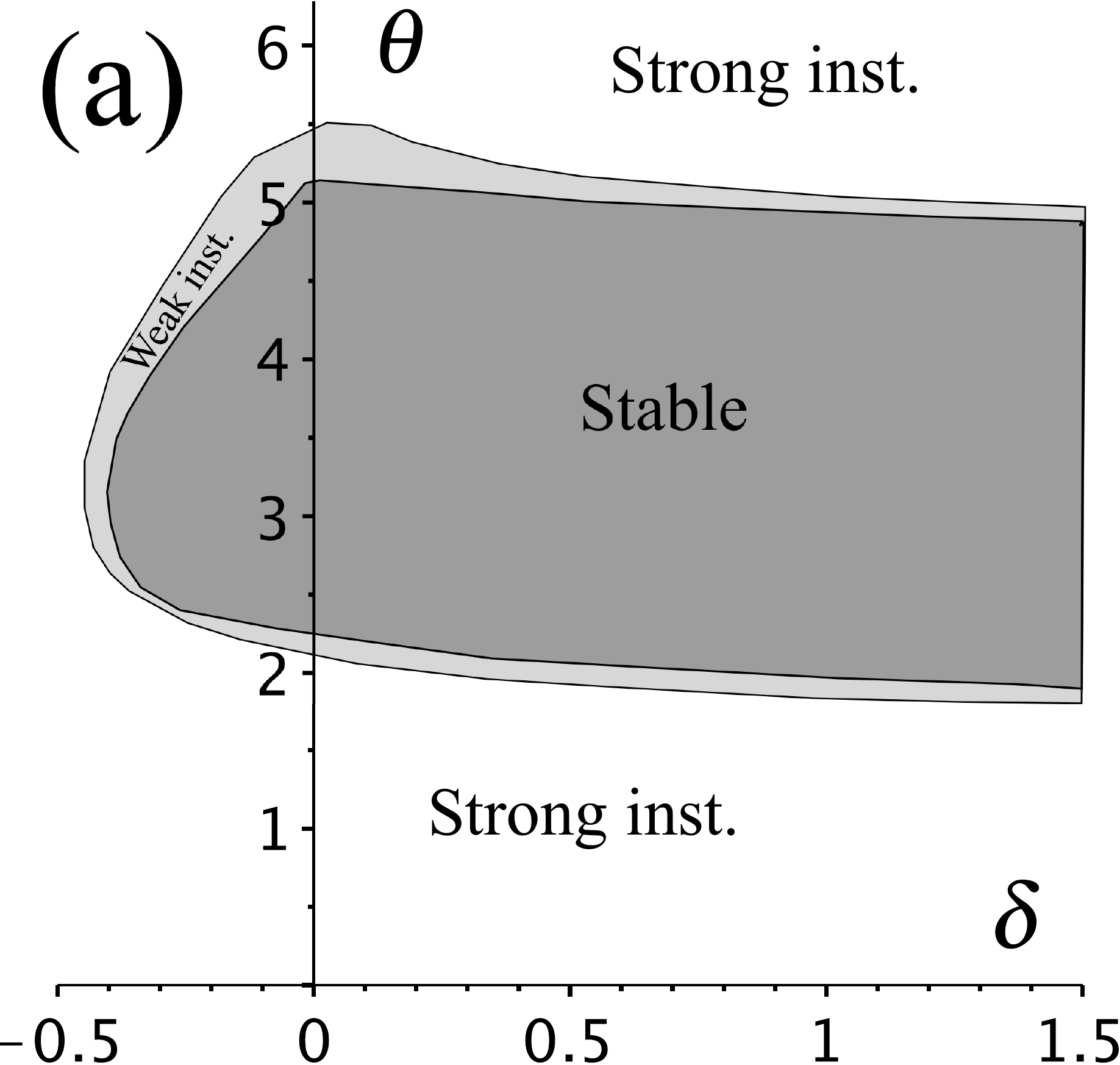}\hskip
1.5 cm\includegraphics[width=0.3\linewidth]{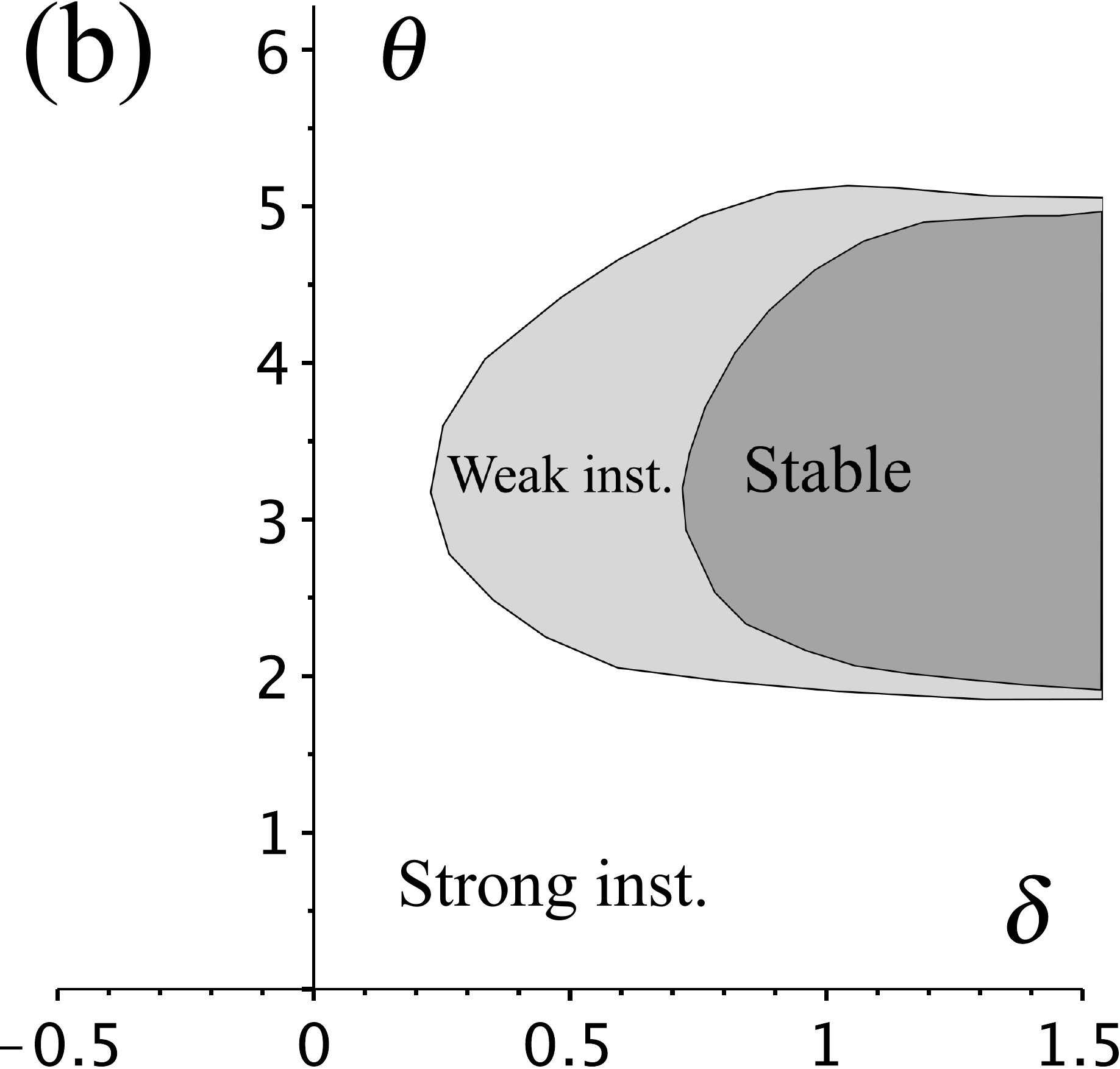}

\caption{Regions of stability (dark gray), weak (light gray) and strong (white)
instability for the plane waves in quintic CGLE on $(\delta,\theta)$
plane. (a) $q=0,$ (b) $q=1$. Other parameters: $\eta=0.2,$ $\beta=0.5,$
$\epsilon=1,$ $\mu=-1$, and $\nu=-0.1$.\label{fig:pwd-bif-wk}}
\end{figure}

Figure \eqref{fig:pwd-sol-wk} shows the upper branch of $\gamma_{1,2}(k,\xi)$
calculated for different values of $\theta$, $\delta$, and $q$.
Even though the surfaces $\gamma(k,\xi)$ are given by explicit expressions,
we were not able to find simple analytical conditions for weak instability
in terms of $\theta$ and $q$. Instead, we determined the sign of
$\sup[\gamma(k,\xi)]$ numerically. The resulting bifurcation diagrams
showing the regions of stable, weakly unstable, and strongly unstable
plane waves with $q=0$ and $q=1$ on $(\delta,\theta)$ plane are
presented in Figure \eqref{fig:pwd-bif-wk}.

\begin{figure}[h]
\centering \includegraphics[width=0.3\linewidth]{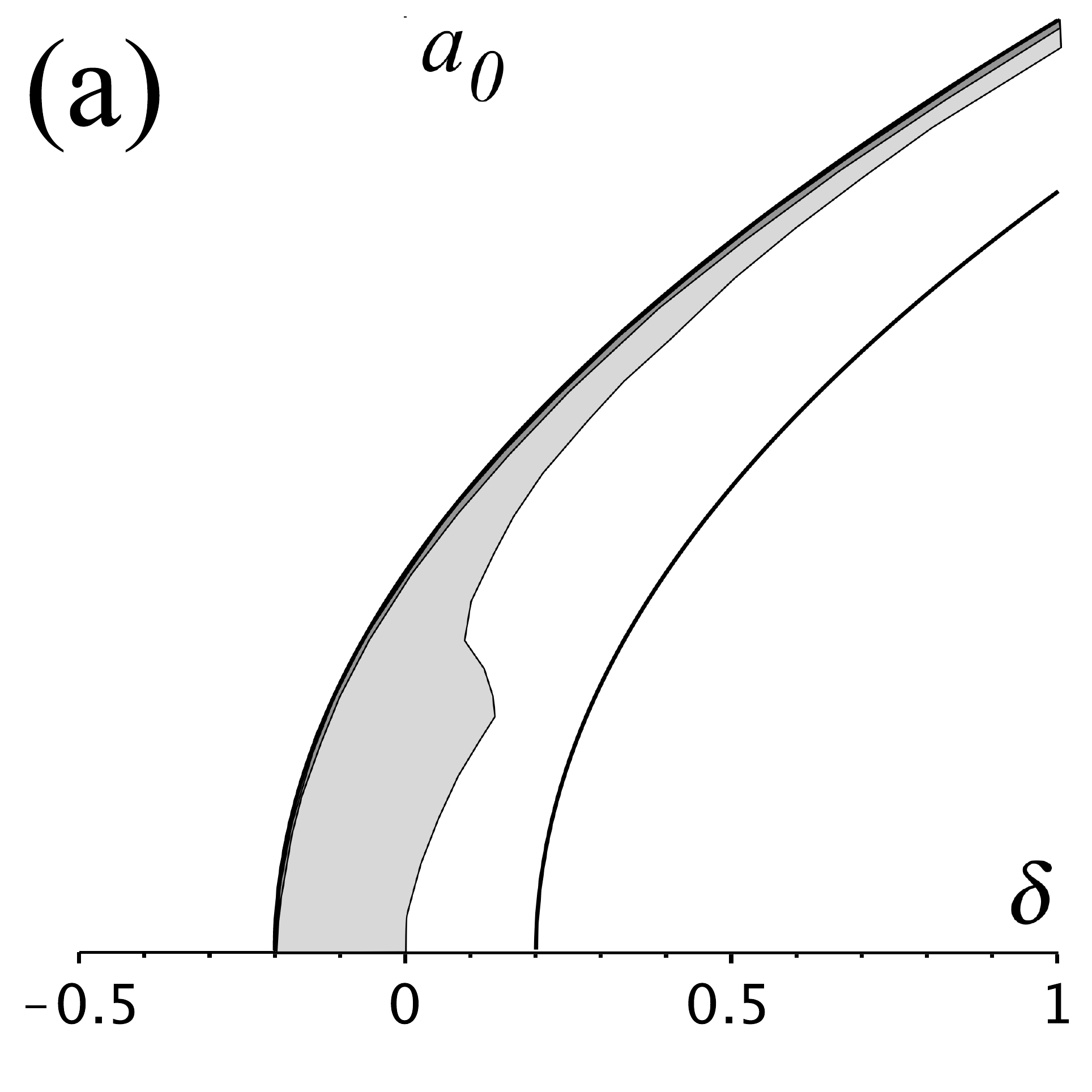}\hskip
1.5 cm\includegraphics[width=0.3\linewidth]{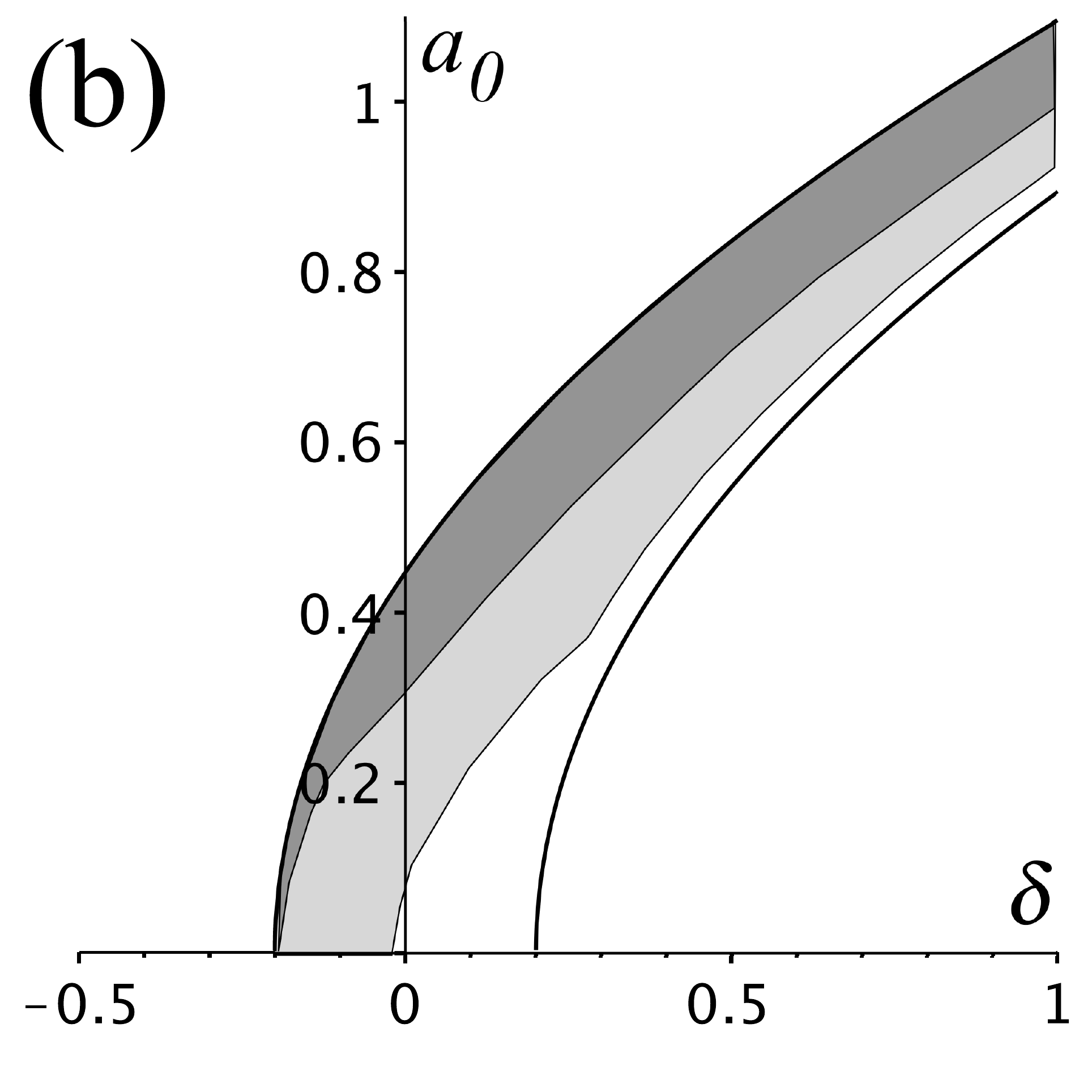}\caption{Regions of stability (dark gray), weak (light gray) and strong (white)
instability for the plane waves in cubic CGLE in $(\delta,a_{0})$
plane. (a) $\theta<\pi$, (b) $\theta>\pi$. Other parameters $q=0,\:\eta=0.2,\:\beta=0.5,$
and $\epsilon=-1$ .\label{fig:Ekh3-pwd-cubic} }
\end{figure}

\begin{figure}[h]
\centering \includegraphics[width=0.3\linewidth]{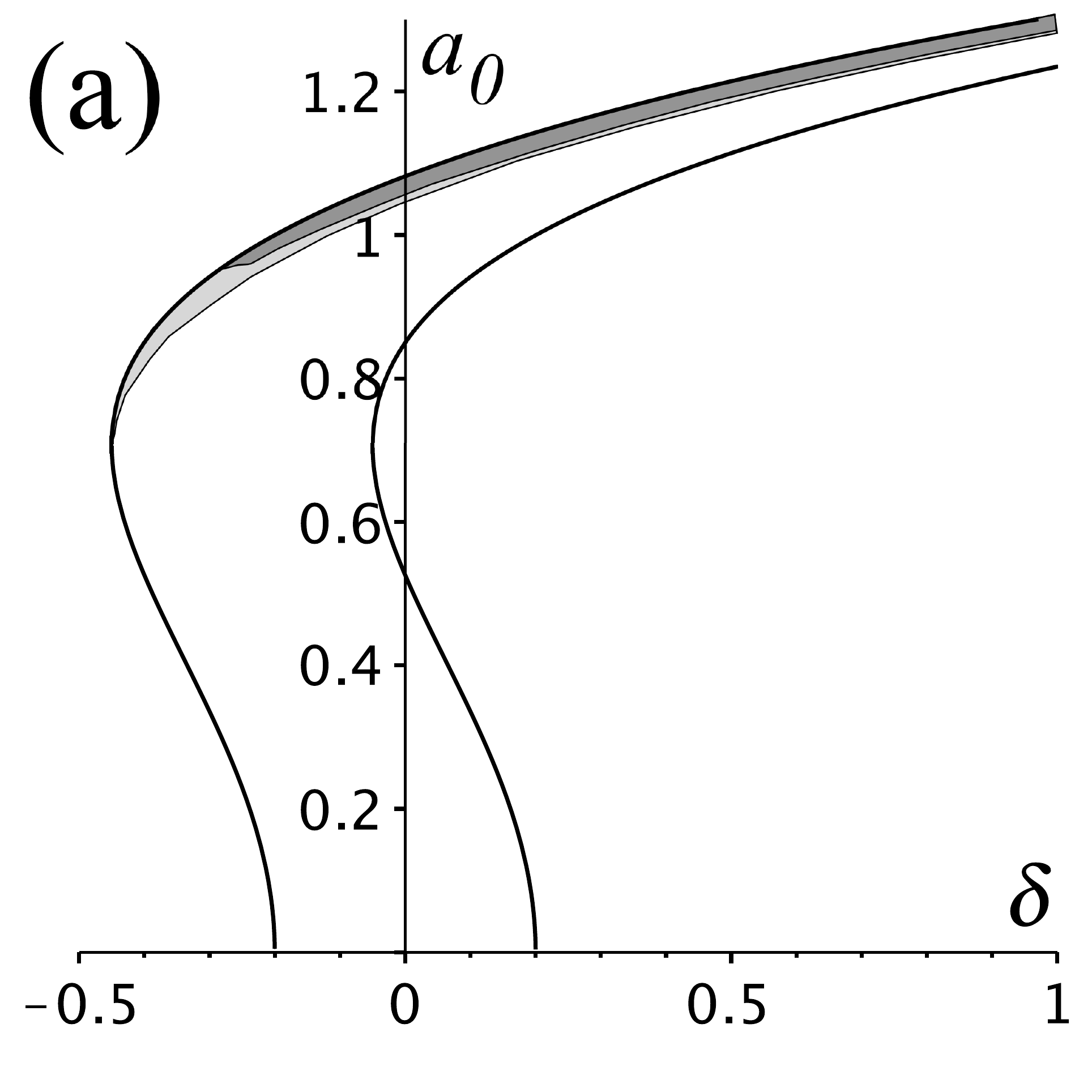}\hskip
1.5 cm\includegraphics[width=0.3\linewidth]{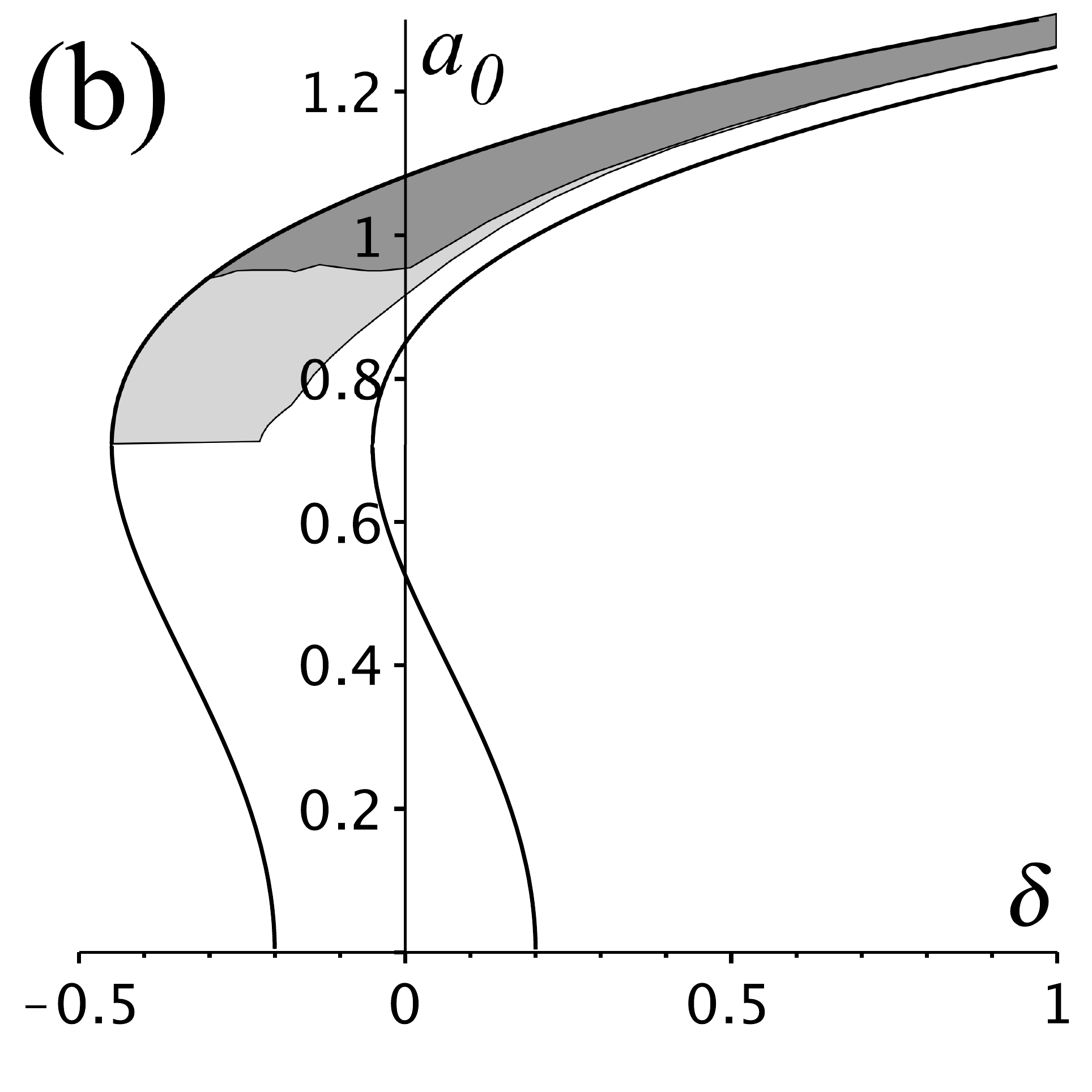}\caption{Regions of stability (dark gray), weak (light gray) and strong (white)
instability for the plane waves in quintic CGLE in $(\delta,a_{0})$
plane. (a) $\theta<\pi$, (b) $\theta>\pi$. Other parameters $q=0,\:\eta=0.2,\:\beta=0.5,$
$\epsilon=1,$ $\mu=-1$, and $\nu=-0.1$.\label{fig:Ekh3-pwd-quintic} }
\end{figure}

The corresponding bifurcation diagrams in $(\delta,a_{0})$ plane
are shown in Figs. \ref{fig:Ekh3-pwd-cubic} and \ref{fig:Ekh3-pwd-quintic}.
Note that in contrast to the unique parameterization of the plane
waves with the parameter $\theta$, in general the plane wave solution
is not uniquely defined by the parameter $a_{0}$ (the frequencies
$\omega$ can still be different for the same $a_{0}$), which means
that the sets with $0\le\theta\le\pi$ and $\pi\le\theta\le2\pi$
are overlapping on the $(\delta,a_{0})$ plane. Moreover, the stability
properties of these two sets are not symmetric with respect to $\theta=\pi$,
as shown in Figs. \ref{fig:Ekh3-pwd-cubic},\ref{fig:Ekh3-pwd-quintic}(a)
and (b), respectively.

Note that the areas of weak instability for quintic CGLE in the large
delay limit are consistent with the stability borders obtained numerically
for $\tau=50$ (compare the transitions from stability to instability
on the branches of solutions and black dashed lines determining the
boundaries in the large delay limit in Fig. \ref{fig:pwd-stab-q0}(b)).
It provides also a simple qualitative way how to predict the appearance
of stable plane waves: namely, when the branches of plane waves for
a finite $\tau$ appear to be in the domain of stability given in
Figs.~\ref{fig:Ekh3-pwd-cubic},\ref{fig:Ekh3-pwd-quintic} (dark
grey shaded), then they are likely to be stable. Since the domains
are independent on $\tau$, the number of coexisting stable plane
waves grows linearly with $\tau$ \cite{PhysRevE.79.046221}.

\subsection{Numerical simulations \label{sub:Numerical-simulations-of}}

To investigate the behavior of particular plane wave solutions, e.g.
the onset of destabilization and convergence to stable solutions,
we performed direct numerical integration of the delayed quintic CGLE
\eqref{CGLE}. An embedded adaptive Cash-Carp scheme for time-stepping
was used with the spatial derivative was treated by three-point central
finite difference scheme. Periodic boundary conditions were are applied,
and the length of the system is chosen to include 16 spatial periods
of simulated plane waves with $q=1$. Space was discretized into $500$
points, while relative tolerance for Cash-Karp scheme was set at $10^{-6}$.
\begin{figure}[h]
\centering \includegraphics[width=0.3\linewidth]{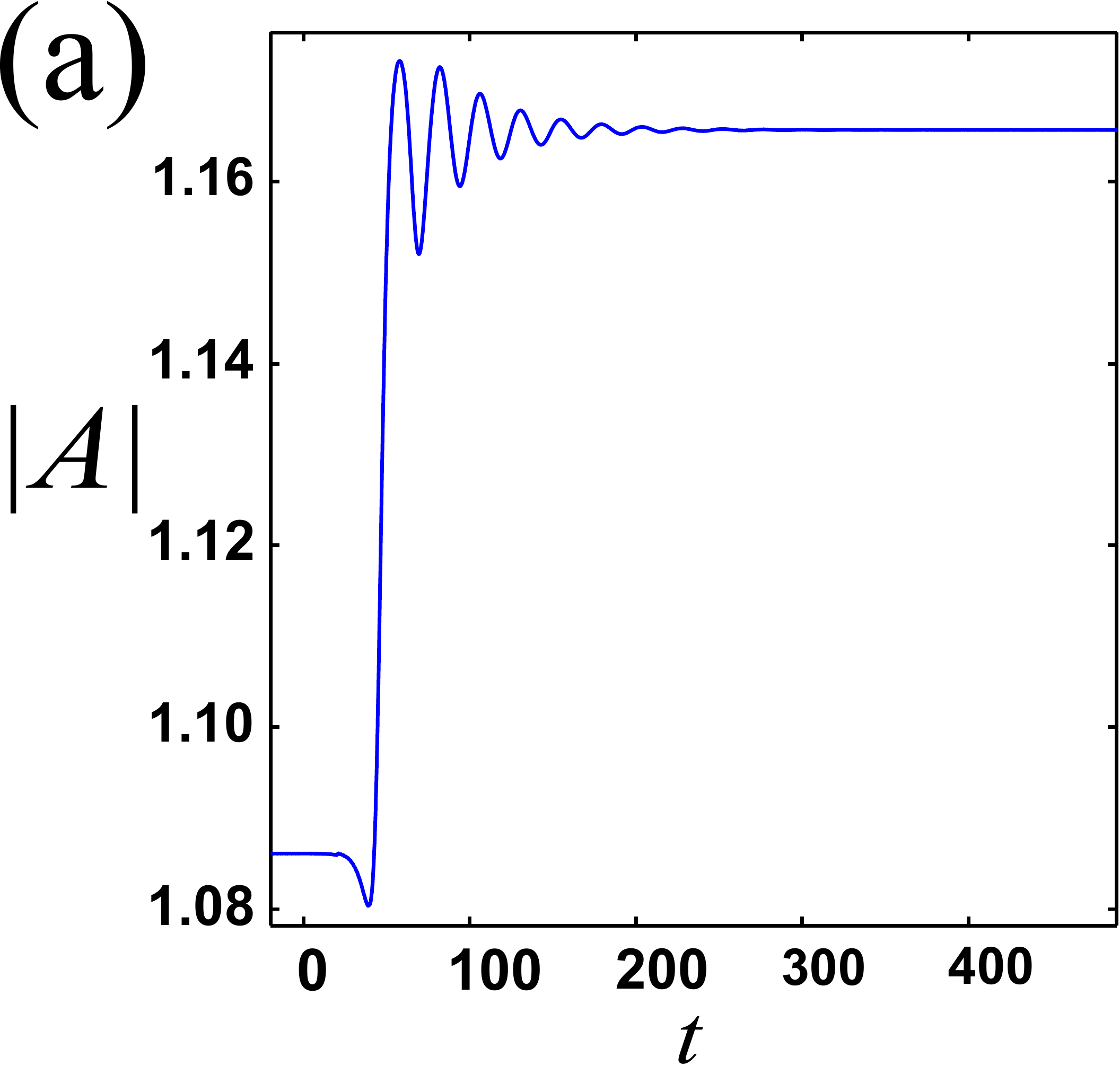}\hskip
1.5 cm\includegraphics[width=0.3\linewidth]{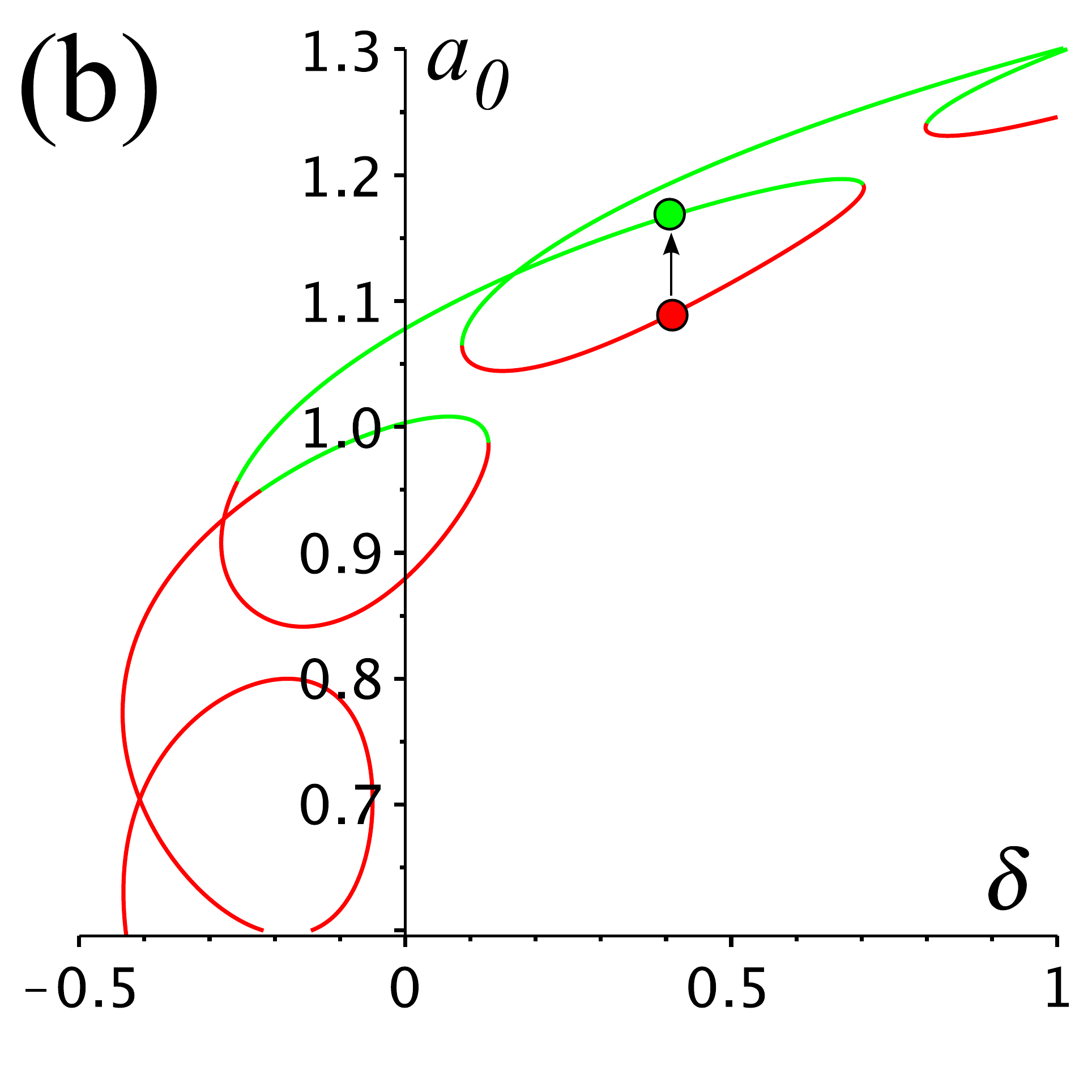}

\caption{Solution starting from the neighborhood of an unstable plane wave
with $q=0$ is attracted to a stable plane wave with the same wavenumber.
(a) depicts the time evolution of the amplitude $a$. (b) shows the
corresponding transition in $(\delta,a_{0})$ plane. Initial unstable
solution: $\omega=1.12,\: a_{0}=1.086$ (red circle). Resulting stable
solution: $\omega=0.998,\: a_{0}=1.166$ (green circle). System parameters:
$\delta=0.023$, $\tau=20,$ $\eta=0.2,$ $\beta=0.5,$ $\epsilon=1,$
$\mu=-1$, and $\nu=-0.1$.\label{fig:pwd-evol_q0}}
\end{figure}

First, we considered the solution starting in the vicinity of the
unstable homogeneous plane wave with $q=0$, see red point in Fig.
\ref{fig:pwd-evol_q0}(b). Fig. \ref{fig:pwd-evol_q0}(a) shows the
evolution of the amplitude, while Fig. \ref{fig:pwd-evol_q0}(b) shows
the stable plane wave (green circle), to which the solution is attracted.
In the case when there are no stable plane waves for given system
parameters, the solution converges to the stable homogeneous state
$A=0$.

\begin{figure}[ph]
\centering \includegraphics[width=0.6\linewidth]{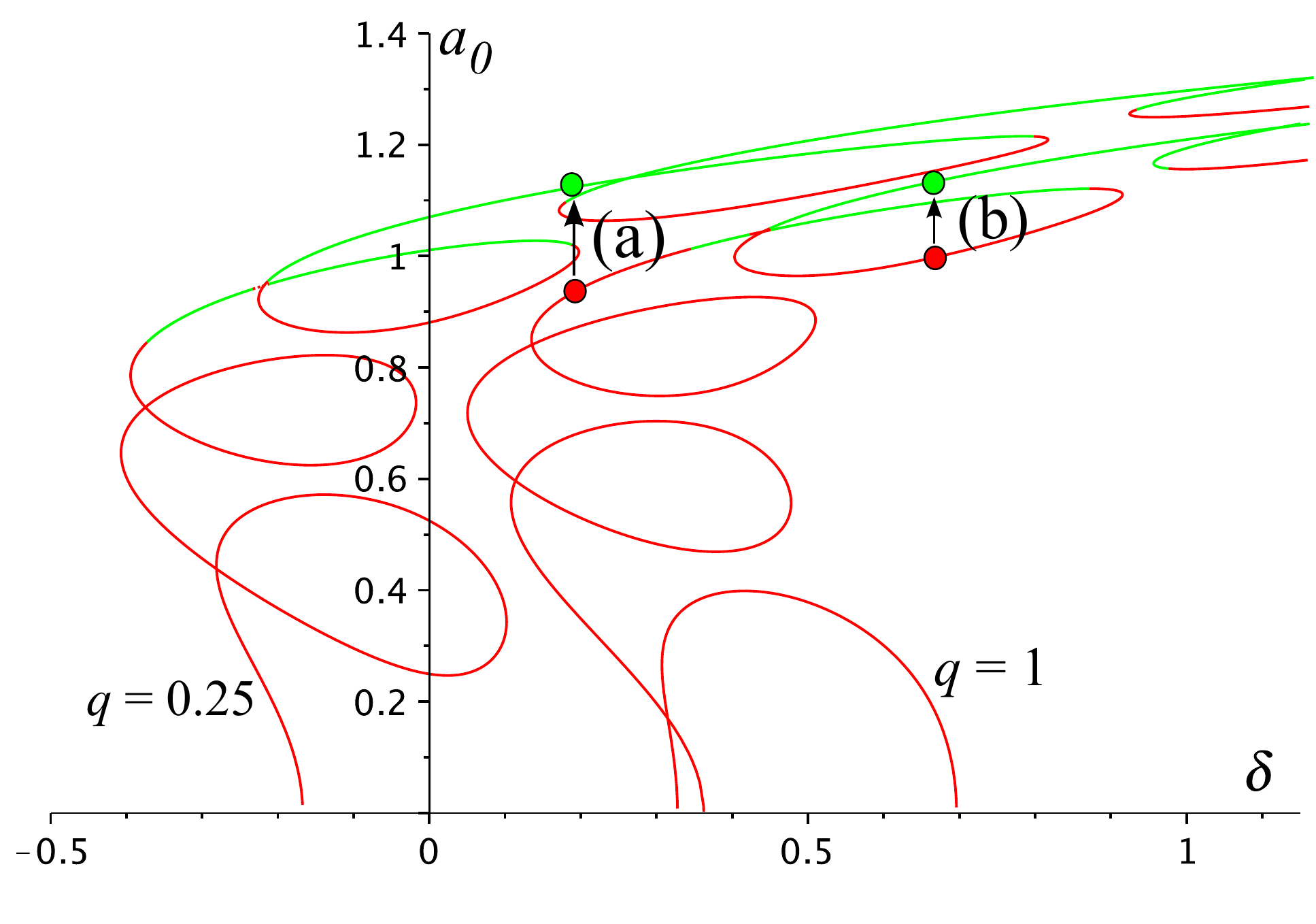}\caption{Evolution of solutions starting close to unstable plane waves with
$q=1$. Solution (a) with $\delta=0.19$, $\omega=0.31,$ and $a_{0}=0.87$
approach the stable plane wave with different wavenumber $q=0.25,$
and $\omega=0.969,$ $a_{0}=1.12$. Solution (b) $\delta=0.66$, $\omega=0.5,$
$a_{0}=0.99$ approaches the stable solution with the same wavenumber
$q=1,$ and $\omega=0.626,$ $a_{0}=1.11$. \label{fig:pwd_evol_q1amp} }
\end{figure}

\begin{figure}[ph]
\centering \includegraphics[width=0.4\linewidth]{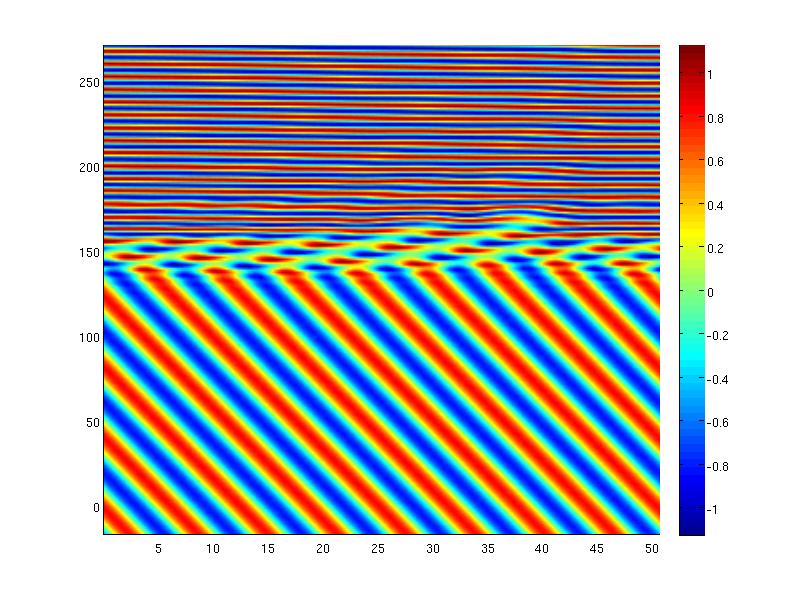}\hskip
1.5 cm\includegraphics[width=0.4\linewidth]{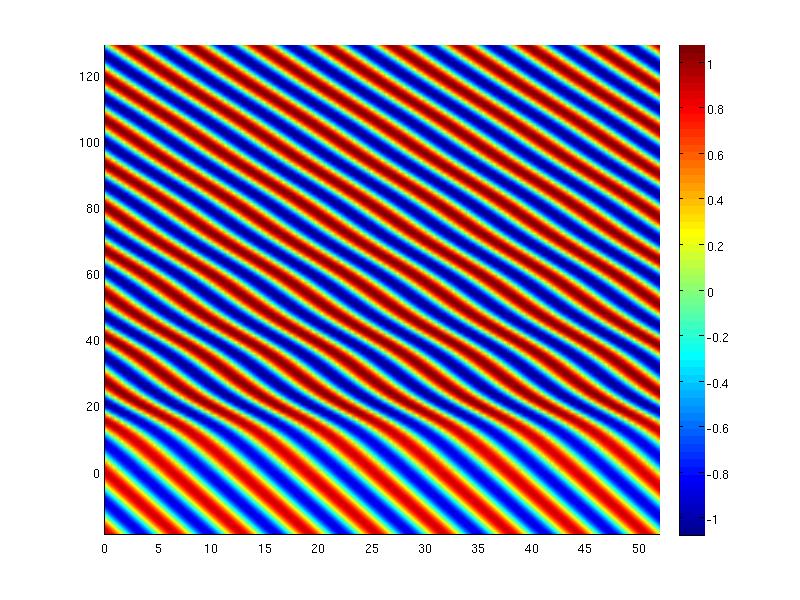}

\caption{Evolution of plane waves with $q=1$. The profile of $\Re[A]$ is
shown in color. System parameters: $\tau=20,$ $\eta=0.2,$ $\beta=0.5,$
$\epsilon=1,$ $\mu=-1$, and $\nu=-0.1$. Panel (a) corresponds to
the transition (a) in Fig.~\ref{fig:pwd_evol_q1amp} when a defect
is developed and the spatial wavenumber is changed. Panel (b) corresponds
to the wavenumber-preserving transition (b) in Fig.~\ref{fig:pwd_evol_q1amp}
without defects. \label{fig:pwd_evol_q1_prof} }
\end{figure}

Transitions similar to those shown in Fig. \ref{fig:pwd-evol_q0}
were also observed for $q\ne0$, see Fig. \ref{fig:pwd_evol_q1amp},
where two unstable plane waves with $q=1$ were chosen as initial
conditions. Figures \ref{fig:pwd_evol_q1_prof}(a) and (b) present
the spatio-temporal evolution of the solutions corresponding to the
transitions shown in Fig. \ref{fig:pwd_evol_q1amp}. It is seen that
the solution (a) develops defects after several delay cycles and transforms
into a slightly modulated solution with lower spatial wavenumber $(q=0.25)$,
which is stable for the given control parameter values. Note that
for these parameter values there are no stable plane waves with the
wavenumber $q=1$. By contrast, solution (b) does not change its principal
spatial wavenumber. After several delay periods of transient it converges
without defects to a stable plane wave having the same wavenumber
$q=1$, see corresponding transition (b) in Fig.~\ref{fig:pwd_evol_q1amp}.

\begin{figure}[h]
\centering \includegraphics[height=0.2\textheight]{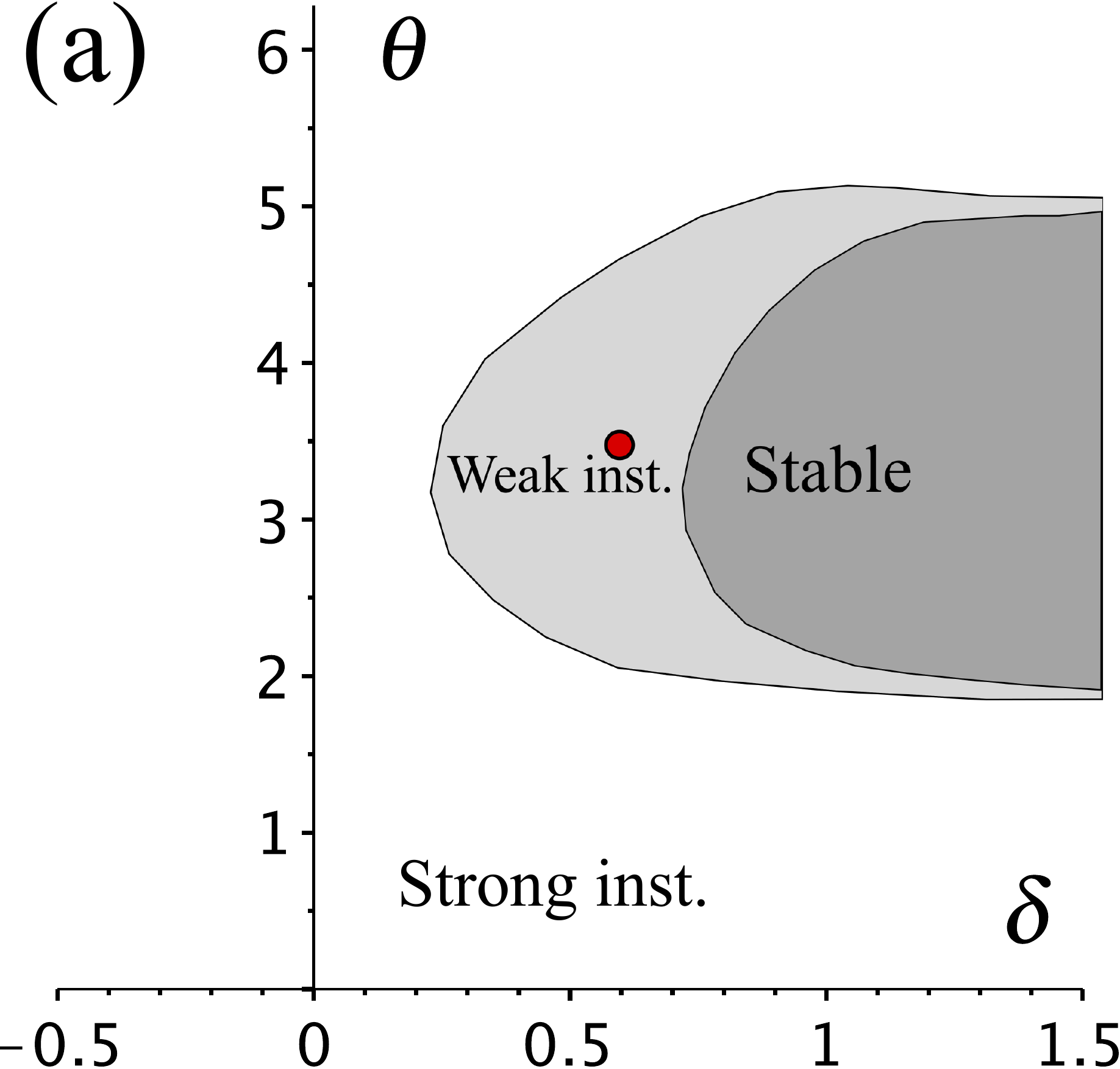}\hskip
1.5 cm\includegraphics[height=0.2\textheight]{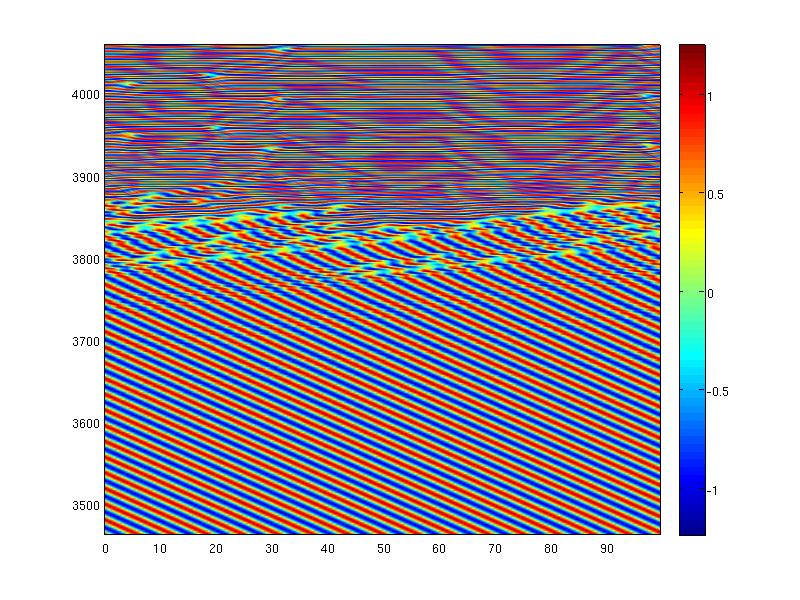}

\caption{A solution started from the vicinity of a weakly unstable plane wave
with $q=1$ for large delay. The weak instability of the plane wave
implies that the perturbation stays small at least several delay intervals.
Initial unstable solution: $\theta=3.94,\:\omega=0.33,\: a_{0}=1.05$
(red circle in (a)). System parameters: $\delta=0.56$, $\tau=50,$
$\eta=0.2,$ $\beta=0.5,$ $\epsilon=1,$ $\mu=-1$, and $\nu=-0.1$.\label{fig:pwd-evol_q1-weak}}
\end{figure}

Now we choose $\tau$ sufficiently large ($\tau=50$), so that we
can exploit the results of asymptotic analysis obtained in the limit
of large delay. We consider a plane wave with the control parameters
from weakly unstable region on $(\theta,\delta)$ plane and investigate
the evolution of a small perturbation of this plane wave. We expect
to observe the onset of destabilization after a long period. Figure
\ref{fig:pwd-evol_q1-weak}(a) illustrates the choice of the initial
plane wave, while Fig. \ref{fig:pwd-evol_q1-weak}(b) shows the transient
and the onset of destabilization. The solution stays close to the
weakly unstable plane wave for about $80$ delay time periods, and
then goes away from it. Eventually, the solution is attracted to a
stable plane wave with $q=1/8$.

\section{Conclusions}

We have investigated the properties of plane wave solutions of cubic
and cubic-quintic CGLE with delayed feedback. It is demonstrated that
the delayed feedback induces a multistability of plane wave solutions
with the same wavenumber $q$. As the gain parameter $\delta$
is varied, the branches of plane wave solutions are shown to exhibit
a snaking behavior, where the frequency of the snaking oscillations
is proportional to the time delay $\tau$. Furthermore, stability
properties of trivial homogeneous zero solution and plane wave solutions
with different wavenumbers are investigated. The numerical bifurcation
diagrams for various delay times as well as the analytical results
in the limit of large delay reveal the borders of strong and weak
instability of the plane waves. Direct numerical integration of the
model equation confirms the results of analytical investigations.

\bibliographystyle{plain}
\bibliography{../delay}

\end{document}